\documentclass{amsart}
\usepackage{latexsym}
\usepackage{amsmath}
\usepackage{amssymb}
\usepackage{amsthm}
\usepackage{amscd}
\usepackage{color}
\usepackage[notref,notcite]{showkeys} 
\numberwithin{equation}{section}
\numberwithin{equation}{section}
\usepackage{fancyhdr}
\definecolor{dblue}{rgb}{0,0,0.45}
\definecolor{red}{rgb}{0.7,0,0}

\if0
\documentclass{amsart}
\usepackage{latexsym,amssymb,amsmath,amsthm,amscd,graphicx}

\fi

\newtheorem{theorem}{Theorem}[section]
\newtheorem{lemma}[theorem]{Lemma}
\newtheorem{corollary}[theorem]{Corollary}
\newtheorem{proposition}[theorem]{Proposition}

\newtheorem{problem}[theorem]{Problem}

\theoremstyle{definition}

\newtheorem{remark}[theorem]{Remark}
\newtheorem{definition}[theorem]{Definition}

\theoremstyle{remark}

\newcommand{\N}{{\mathbb N}}
\newcommand{\R}{{\mathbb R}}
\newcommand{\C}{{\mathbb C}}
\newcommand{\Q}{{\mathbb Q}}

\begin{document}

\title{The singular integral operator and its commutator on weighted Morrey spaces}

\author{Shohei Nakamura and Yoshihiro Sawano}

\maketitle

\begin{abstract}
In this paper, 
we give necessary conditions and sufficient conditions respectively for 
the boundedness of the singular integral operator 
on the weighted Morrey spaces. 
We observe the phenomenon unique to the case of Morrey spaces;
the $A_q$-theory by Muckenhoupt and Wheeden does not suffice. 
We also discuss the boundedness of the commutators. 
The difference between the maximal operator and the singular integral operators
will be clarified.
Further, the property of the sharp maximal operator
is investigated.
\end{abstract}

\section{Introduction}
\label{s1}
The weight theory on Lebesgue spaces have been developed 
by many mathematicians in terms of $A_p$ conditions 
and is also developing; see for example \cite{CF74, Hytonen11-0, MH76}. 
As a remarkable recent progress, 
Hyt\"{o}nen \cite{Hytonen11-0} solved the ``$A_2$ conjecture" 
which was an open problem on this field.  
Meanwhile, some mathematicians pay attention to weighted Morrey spaces and 
try to construct the weight theory on Morrey spaces; for example \cite{KoSh09, NS15, NS09, Ta15}. 
As is pointed out in \cite{NS09, Ta15} implicitly,  
$A_p$ condition is not suitable on weighted Morrey spaces and hence, 
it is a natural problem to find the alternative condition. 
To tackel this problem, 
we discuss the boundedness of the singular integral operators
on weighted Morrey spaces 
by assuming the boundedness of the maximal operator in this paper.  

First of all, we fix our notation.
We denote the set of the all Lebesgue measurable functions by $L^0(\R^n)$.
By $\mathcal{Q}=\mathcal{Q}(\R^n)$, 
we mean the all cubes in $\R^n$ whose sides parallel to the
coordinate axes.
We denote the family of all dyadic cubes by $\mathcal{D}=\mathcal{D}(\R^n)$ 
and the family of all dyadic cubes with respect to $Q\in\mathcal{Q}$ 
by $\mathcal{D}(Q)$.
A weight $w$ is a locally integrable function on $\R^n$ such that 
$w(x)>0$ for almost everywhere $x\in\R^n$.
For weight $w$ and measurable set $E\subset\R^n$, 
we denote 
$w(E):=\int_Ew(x)dx$.
In addition, we denote a weighted measure by $dw$, 
that is, $dw(x)=w(x)dx$.

We define the weighted Morrey space 
$\mathcal{M}^p_q(w_1,w_2)$, 
where $0< q\leq p<\infty$ and $w_1$, $w_2$ are weights.
The weighted Morrey space 
$\mathcal{M}^p_q(w_1,w_2)$ is the set of all functions
$f\in L^q_{\rm loc}(w_2)$ for which the quasi-norm
\begin{equation*}
\|f\|_{\mathcal{M}^p_q(w_1,w_2)}
:=
\sup_{Q\in\mathcal{Q}}w_1(Q)^{\frac{1}{p}-\frac{1}{q}}
\left(\int_{Q}|f(x)|^qdw_2(x)\right)^\frac{1}{q}
\end{equation*}
is finite.
When $w_1=w_2=w$, 
$\mathcal{M}^p_q(w_1,w_2)$ corresponds to 
the weighted Morrey space $\mathcal{M}^p_q(w,w)$ 
introduced by Komori and Shirai 
in \cite{KoSh09}.
Meanwhile, 
when $w_1=dx$ and $w_2=w$, 
$\mathcal{M}^p_q(w_1,w_2)$ corresponds to 
the weighted Morrey space $\mathcal{M}^p_q(dx,w)$ 
introduced by Samko 
in \cite{NS09}. 
Our main interest is the latter case: $\mathcal{M}^p_q(dx,w)$. 

In \cite{NS15}, 
we introduced the weight class $\mathcal{B}_{p,q}$ 
and the weighted integral condition
in the context of the boundedness of 
the Hardy-Littlewood maximal operator 
on weighted Morrey spaces of Samko type.
\begin{definition}[\cite{NS15}]
Let $0<q\leq p<\infty$ and $w$ be a weight.
\begin{enumerate}
\item
One says that a weight $w$ is in the class $\mathcal{B}_{p,q}$ 
if there exists $C_{p,q}>0$ such that for any $Q_0\in\mathcal{Q}$,
\begin{equation}\label{150422-1}
\sup_{Q\in\mathcal{Q}:Q\subset Q_0}
\Phi_{p,q,w}(Q)\leq
C_{p,q}\Phi_{p,q,w}(Q_0),
\end{equation} 
or equivalently, 
\begin{equation}\label{160627-1}
\|\chi_{Q_0}\|_{\mathcal{M}^p_q(dx,w)}
\sim
\Phi_{p,q,w}(Q_0)
\end{equation}
hold,
where we defined 
\begin{equation*}
\Phi_{p,q,w}(Q):=
|Q|^\frac{1}{p}\left(\frac{w(Q)}{|Q|}\right)^\frac{1}{q}
\quad (Q \in {\mathcal Q}).
\end{equation*}
\item
The weighted integral condition for $p,q$ and $w$ holds,
if there exists a constant $C>0$ such that
\begin{equation}\label{weight-integral}
\int_1^\infty\frac{1}{\Phi_{p,q,w}(sQ)}\frac{ds}{s}
\leq
\frac{C}{\Phi_{p,q,w}(Q)}
\quad
(Q\in\mathcal{Q})
\end{equation}
holds.
\end{enumerate}
\end{definition}
Note that the weighted integral condition (\ref{weight-integral})
for $p,q$ and $w$ implies $w\in \mathcal{B}_{p,q}$; 
see \cite{NS15} for details.
Using these two notions, 
we aim to give necessary conditions and sufficiently conditions 
for the boundedness of 
the singular integral operators.
To this end, we first recall the definition of 
the Hardy-Littlewood maximal operator and the K\"{o}the dual space 
of Morrey spaces.
By $M$, we mean the (unweighted) Hardy-Littlewood maximal operator:
\begin{equation*}
Mf(x):=\sup_{Q\in\mathcal{D}}\frac{1}{|Q|}\int_Q |f(y)|dy\cdot\chi_{Q}(x).
\end{equation*}
To recall the K\"{o}the dual spaces of Morrey spaces, 
we need some terminologies. 
Let $0<\alpha\leq n$ and $E\subset\R^n$
be an arbitrary set.
Then the $\alpha$-dimensional Hausdorff content of $E$ 
is defined by
\begin{equation*}
H^\alpha(E)
:=
\inf
\left\{\sum_j\ell(Q_j)^\alpha
:\{Q_j\}_{j}\subset\mathcal{Q}, 
\ 
\bigcup_{j}Q_j\supset E\right\}.
\end{equation*}
For a non-negative function $\phi\geq0$, 
the Choquet integral of $\phi$ with respect to the Hausdorff content
$H^\alpha$ is defined by
\begin{equation*}
\int_{\R^n}
\phi dH^\alpha
:=
\int_0^\infty
H^\alpha
(\{x\in\R^n:\phi(x)>t\})
dt.
\end{equation*}
\begin{definition}[\cite{AX12}]
Define  
$
\mathfrak{B}_\alpha
:=
\left\{
b\in A_1:\int_{\R^n}b dH^\alpha
\leq1
\right\}
$
for $0<\alpha<n$.
\end{definition}
A typical example of $\mathfrak{B}_\alpha$ 
is 
\begin{equation}\label{160111-1}
b_Q
:=
\frac{(M\chi_Q)^{\alpha/n+\varepsilon}}{\ell(Q)^\alpha},
\end{equation}
where $Q\in\mathcal{Q}$ and $\varepsilon \in (0,1-\frac{\alpha}{n})$ 
are arbitrary.
Then for $1<q\leq p<\infty$, 
the space $H^{q',n(1-q/p)}(dx,dx)=H^{q',n(1-q/p)}$ is defined by 
the set of all measurable functions $f$ for which norm 
\begin{equation*}
\|f\|_{H^{q',n(1-q/p)}}
:=
\inf_{b\in \mathfrak{B}_{n(1-q/p)}}
\left(\int_{\R^n}|f(x)|^{q'}b(x)^{-\frac{q'}{q}}dx
\right)^\frac{1}{q'}
\end{equation*}
is finite.
\if0
Since we will use
the unweighted type space $H^{q',n(1-q/p)}(dx,dx)$, 
we abbreviate it by $H^{q',n(1-q/p)}$ in this paper.
\fi

The following are proved in \cite{AX12} or \cite[(2.5)]{Ta15}:
For $1<q\leq p<\infty$, we have 
\begin{equation}\label{dual}
\left|\int_{\R^n}f(x)g(x)dx\right|
\leq
C\|f\|_{\mathcal{M}^p_q(dx,dx)}
\|g\|_{H^{q',n(1-q/p)}},
\end{equation}
for all $f \in \mathcal{M}^p_q(dx,dx)$
and
$g \in H^{q',n(1-q/p)}$
and that
\begin{equation}\label{dual2}
\|g\|_{H^{q',n(1-q/p)}}
\sim
\sup\left\{
\|f \cdot g\|_{L^1}
\,:\,
f\in\mathcal{M}^p_q(dx,dx):\\ \|f\|_{\mathcal{M}^p_q(dx,dx)}
\leq 1\right\},
\end{equation}
for any measurable function $g$.

Now, we state our main results in this paper. 
The Riesz transform is defined by 
\begin{equation*}
R_if(x)
:=
\lim_{\varepsilon\to0}
\int_{\left\{y\in\R^n:|x_i-y_i|>\varepsilon\right\}}
\frac{x_i-y_i}{|x-y|^{n+1}}f(y)dy
\quad
(f\in L^2(\R^n)),
\end{equation*}
for $i=1,\ldots, n$.

\begin{theorem}\label{th-160626-1}
Let $1<q\leq p<\infty$ and $w$ be a weight. 
\begin{enumerate}
\item\label{item-160626-1}
(Sufficient condition)
Assume that the Hardy-Littlewood maximal operator $M$ is bounded on 
$\mathcal{M}^p_q(dx,w)$ and the weighted integral condition \eqref{weight-integral} 
for $p,q$ and $w$. 
Then for any $f\in\mathcal{M}^p_q(dx,w)$, 
$R_if$ is well defined and $R_i$ is bounded on $\mathcal{M}^p_q(dx,w)$: 
\[
\left\|R_if\right\|_{\mathcal{M}^p_q(dx,w)}
\leq 
C
\|f\|_{\mathcal{M}^p_q(dx,w)}
\quad
(f\in\mathcal{M}^p_q(dx,w)). 
\]
\item\label{item-160626-2}
(Necessary condition)
Conversely, assume $R_i$ is bounded on $\mathcal{M}^p_q(dx,w)$. 
Then we have 
\begin{equation}\label{151221-1}
\frac{1}{|Q|}
\left\|
\chi_Qw^\frac{1}{q}
\right\|_{\mathcal{M}^p_q(dx,dx)}
\left\|
\chi_{Q}
w^{-\frac{1}{q}}
\right\|_{H^{q',n(1-q/p)}}
\leq
C
\quad
(Q\in\mathcal{Q}),
\end{equation}
and the weighted integral condition \eqref{weight-integral} 
for $p,q$ and $w$.
\end{enumerate} 
\end{theorem}

\begin{remark}
In Theorem \ref{th-160626-1}, 
we have a ``gap" between the sufficient condition \eqref{item-160626-1} and 
the necessary condition \eqref{item-160626-2}, 
that is, a ``gap" between the boundedness of $M$ and \eqref{151221-1}. 
To this problem, 
it was conjectured in \cite{Ta15} that 
the boundedness of $M$ on $\mathcal{M}^p_q(dx,w)$ is equivalent to 
\eqref{151221-1}. 
Hence, if the conjecture is true, then 
we obtain from Theorem \ref{th-160626-1} 
the characterization of the boundedness of $R_i$ on $\mathcal{M}^p_q(dx,w)$. 
That is, $R_i$ is bounded on $\mathcal{M}^p_q(dx,w)$ 
if and only if 
\eqref{151221-1} and the weighted integral condition \eqref{weight-integral} hold. 
\end{remark} 

As is well known, 
on the weighted Lebesgue spaces, 
the boundedness of $R_i$ or $M$ is equivalent to $A_p$ condition. 
However, as in Theorem \ref{th-160626-1}, 
there exists an essential difference between the behavior of $R_i$ and $M$ on the weighted Morrey spaces 
coming from the weighted integral condition \eqref{weight-integral}. 
%Such a phenomenon is unique to the Morrey spaces. 
Of course, we have possibility 
that the weighted integral condition \eqref{weight-integral} is not essential 
which means \eqref{weight-integral} always holds automatically. 
In fact, in the case of weighted Lebesgue space: $p=q$, 
the weighted integral condition \eqref{weight-integral} holds as long as 
$w$ is a doubling weight. 
To see that the weighted integral condition 
\eqref{weight-integral} is essential,
we focus on the specific weight, namely, the power weight: $w(x)=w_\alpha(x):=|x|^{\alpha}$. 
On this specific situation, 
the problem of the boundedness of $M$ on $\mathcal{M}^p_q(dx,w)$ 
is completely solved in \cite{Ta15} as follows. 

\begin{theorem}[\cite{Ta15}]\label{th-tanaka}
Let $1<q< p<\infty$ and $\alpha>-n$. 
Then the following are equivalent:
\begin{enumerate}
\item
The boundedness of $M$ on $\mathcal{M}^p_q(dx,w_\alpha)$: 
\[
\left\|Mf\right\|_{\mathcal{M}^p_q(dx,w_\alpha)}
\leq
C
\|f\|_{\mathcal{M}^p_q(dx,w_\alpha)}
\quad
(f\in\mathcal{M}^p_q(dx,w_\alpha));
\]
\item
The dual inequality \eqref{151221-1} with $w=w_\alpha$: 
\[
\frac{1}{|Q|}
\left\|
\chi_Qw_\alpha{}^\frac{1}{q}
\right\|_{\mathcal{M}^p_q(dx,dx)}
\left\|
\chi_{Q}
w_\alpha{}^{-\frac{1}{q}}
\right\|_{H^{q',n(1-q/p)}}
\leq
C
\quad
(Q\in\mathcal{Q});
\]
\item
The range of $\alpha$: 
\begin{equation}\label{160117-1}
-\frac{q}{p}n
\leq
\alpha
<
n\left(q-\frac{q}{p}\right).
\end{equation}
\end{enumerate}
\end{theorem}

Meanwhile, as is shown in \cite{NS15}, 
the weighted integral condition \eqref{weight-integral} for 
$p,q$ and $w_\alpha$ is equivalent to 
$\alpha>-\frac{q}{p}n$. 
Hence, we obtain the following:
\begin{corollary}\label{cr-160113-2}
Let $1<q\leq p<\infty$ and $\alpha>-n$. 
The Riesz transform $R_i$, $i=1,\ldots,n$ is bounded on $\mathcal{M}^p_q(dx,w_\alpha)$
if and only if 
\begin{equation}\label{160121-5}
-\frac{q}{p}n<\alpha<n\left(q-\frac{q}{p}\right).
\end{equation}
\end{corollary}

Comparing \eqref{160117-1} and \eqref{160121-5}, 
we notice the essential difference appears 
at the end-point: $\alpha=-\frac{q}{p}$.

We now compare our result to the known results. 
In \cite{NS15}, we proved the boundedness of $R_i$ on 
$\mathcal{M}^p_q(dx,w)$ 
by assuming $w\in A_q$ and the weighted integral condition \eqref{weight-integral}. 
In particular, when the case of $w=w_\alpha$, 
the conditions: $w\in A_q$ and \eqref{weight-integral} 
are combined into:
\begin{equation}\label{A_q}
-\frac{q}{p}n<\alpha<n(q-1). 
\end{equation}
Hence, by comparing \eqref{160121-5} and \eqref{A_q}, 
we notice that 
Theorem \ref{th-160626-1} improves the result obtained in \cite{NS15}. 
When the case of $n=1$, in \cite{NS09}, 
Samko obtained the characterization of the 
boundedness of the Hilbert transform with 
power weights.
Recall that the Hilbert transform $H$ is a singular integral 
operator defined by 
\begin{equation*}
Hf(x)
:=
\lim_{\varepsilon\to0}
\int_{\left\{y\in\R:|x-y|>\varepsilon\right\}}
\frac{f(y)}{x-y}dy
\quad
(f\in L^2(\R)).
\end{equation*}
Note that the Hilbert transform $H$ is the special case of 
the Riezs transform $R_i$ with $n=1$. 
We can recapture the result by Samko.
\begin{theorem}[\cite{NS09}]\label{NS}
Let $1<q\leq p<\infty$ and $\alpha>-1$.
The Hilbert transform $H$ is bounded on $\mathcal{M}^p_q(dx,w_\alpha)$, 
if and only if 
\begin{equation}\label{160121-4}
-\frac{q}{p}<\alpha<q-\frac{q}{p}.
\end{equation}
\end{theorem}
We see that 
both Theorem \ref{th-160626-1} and Corollary \ref{cr-160113-2} 
cover Theorem \ref{NS}. 

We also discuss the boundedness of the commutators. 
The commutator of the Riesz transform $R_i$  
with a locally integrable function $b$ 
is initially defined by 
\begin{equation*}
[b,R_i]f:=b\cdot R_if-R_i(b\cdot f),
\end{equation*}
for $f\in L^\infty_{\rm c}$.
A locally integrable function $b$ is said to be a
BMO function, if 
\begin{equation}\label{BMO}
\|b\|_{\rm BMO}
:=\|b^\sharp\|_{L^\infty}<\infty,
\end{equation}
where $b^\sharp$ is the Fefferman-Stein-Stromberg sharp maximal function
defined by (\ref{160627-2}).
When $b$ is a BMO function, 
we can extend the commutator $[b,R_i]$ as a bounded linear operator 
on $L^p(\R^n)$ for all $1<p<\infty$; see \cite[Theorem 3.5.6]{grafakos-book} for example.
We also note that in the endpoint case: $p=1$, 
the commutators are more singular 
than the singular integral operators.
%the situation between the singular integral operator and 
%the commutator is different. 
We refer \cite{Perez95} for the result of the commutators in the endpoint case: $p=1$.
Let
$R_i$, $i=1,\ldots, n$ be the $i$-th Riesz transform. 
For the commutators, we obtain a similar result as follows. 
\begin{theorem}\label{th-160113-1}
Let $1<q\leq p<\infty$, $b\in {\rm BMO}$ and 
$i=1,2,\ldots,n$ be fixed.
Assume that 
$M$ is bounded on $\mathcal{M}^p_q(dx,w)$ and the weighted integral condition
(\ref{weight-integral}) for $p,q$ and $w$. 
Then we can extend the commutator $[b,R_i]$ 
as a bounded linear operator 
on $\mathcal{M}^p_q(dx,w)$:
\begin{equation*}
\|[b,R_i]f\|_{\mathcal{M}^p_q(dx,w)}
\leq
C\|f\|_{\mathcal{M}^p_q(dx,w)}
\quad
(f\in\mathcal{M}^p_q(dx,w)).
\end{equation*}
\end{theorem}

To prove Theorem \ref{th-160626-1}, 
we will employ the sharp maximal operator arguments. 
To this end, 
we generalize the sharp maximal inequality on Lebesgue spaces to 
our Morrey spaces. 
Our results related to the sharp maximal inequality 
improve the well known sharp maximal inequality 
on Lebesgue spaces in some sense. 
We recall the local sharp maximal operator 
introduced in \cite{Jo65,St79}.
\begin{definition}[\cite{Hytonen11, Lerner13-2}]
Let $f\in L^0(\R^n)$ and $Q\in\mathcal{Q}$.
\begin{enumerate}
\item
The decreasing rearrangement of $f$ on $\R^n$
is defined by
\begin{equation*}
f^*(t):=\left\{\rho>0:
\left|\left\{x\in\R^n:|f(x)|>\rho\right\}\right|<t
\right\},
\quad
(0<t<\infty).
\end{equation*} 
\item
The local mean oscillation of $f$ on $Q$ 
is defined by
\begin{equation*}
\omega_\lambda(f;Q):=
\inf_{c\in\C}\left((f-c)\chi_Q\right)^*(\lambda|Q|),
\quad
(\lambda\in(0,2^{-1})).
\end{equation*}
\item
Assume that the function $f$ is real-valued.
The median of $f$ over $Q$ denoted by $m_f(Q)$
is a real number satisfying that
\begin{equation*}
\left|\left\{x\in Q:f(x)>m_f(Q)\right\}\right|,
\quad
\left|\left\{x\in Q:f(x)<m_f(Q)\right\}\right|
\leq
\frac{1}{2}.
\end{equation*}
Note that the median $m_f(Q)$ is possibly non-unique.
\end{enumerate}
\end{definition}

For $\lambda\in(0,2^{-1})$ and $Q_0\in\mathcal{Q}$, 
the dyadic local sharp maximal operator $M^{\sharp,d}_{\lambda;Q_0}$ 
is defined by 
\begin{equation*}
M^{\sharp,d}_{\lambda;Q_0}f(x)
:=
\sup_{Q\in\mathcal{D}(Q_0)}
\omega_\lambda(f;Q)\chi_Q(x),
\quad
(x\in\R^n, f\in L^0(\R^n)).
\end{equation*}
In this paper, we are interested in the following sharp maximal operator:
\begin{equation*}
M^{\sharp,d}_{\lambda}f(x)
:=
\sup_{Q_0\in\mathcal{Q}}
\sup_{Q\in\mathcal{D}(Q_0)}
\omega_\lambda(f;Q)\chi_Q(x)
\quad
(x\in\R^n, f\in L^0(\R^n)).
\end{equation*}
As is well known, the Fefferman-Stein sharp maximal operator 
defined by
\begin{equation*}
f^{\sharp,\eta}(x)
:=
\sup_{Q\in\mathcal{Q}}
\left(\frac{1}{|Q|}\int_Q|f(y)-f_Q|^\eta dy\right)^\frac{1}{\eta}
\chi_Q(x),
\end{equation*}
where $\eta>0$ and $f_Q:=\frac{1}{|Q|}\int_Qf(x)dx$, 
is a useful tool to estimate the oscillation of the function $f$.
When $\eta=1$, we abbreviate $f^{\sharp,1}$ to $f^{\sharp}$.
\begin{equation}\label{160627-2}
f^{\sharp}(x)
:=
\sup_{Q\in\mathcal{Q}}
\frac{1}{|Q|}\int_Q|f(y)-f_Q| dy\times
\chi_Q(x),
\end{equation}
There exists a close relation between these two types of the sharp maximal operators 
provided by Jawerth and Torchinsky in \cite{Jawa85}:
\begin{equation}\label{151201-4}
M^{(\eta)}M^{\sharp,d}_{\lambda}f(x)
\sim_{\eta,\lambda}
f^{\sharp,\eta}(x),
\end{equation}
for sufficiently small $\lambda$, 
where $M^{(\eta)}$ 
denotes the powered Hardy-Littlewood maximal operator
defined by 
\begin{equation*}
M^{(\eta)}f(x)
:=
\sup_{Q\in\mathcal{Q}}
\left(\frac{1}{|Q|}\int_Q|f(x)|^\eta dx\right)^\frac{1}{\eta}\chi_Q(x).
\end{equation*}
When $\eta=1$, 
$M^{(1)}$ is the original Hardy-Littlewood maximal operator.

In \cite{Lerner13-2}, Lerner proved the following theorem:
\begin{theorem}[\cite{Lerner13-2}]\label{th-151129-1}
Let $f\in L^0(\R^n)$ and $Q_0\in\mathcal{Q}$.
Then there exists a sparse family of 
$\{Q^k_j\}_{k\in\N_0,j\in J_k}\subset\mathcal{D}(Q_0)$
such that for a.e. $x\in Q_0$,
\begin{equation*}
|f(x)-m_f(Q_0)|
\leq
4M^{\sharp,d}_{\lambda_n;Q_0}f(x)
+
2\sum_{k\in\N_0}\sum_{j\in J_k}\omega_{\lambda_n}(f;Q^k_j)\chi_{Q^k_j}(x).
\end{equation*}
Here, $\lambda_n:=2^{-n-2}$.
\end{theorem}
Here, we say that 
the family $\{Q^k_j\}_{k\in\N_0,j\in J_k}$ 
is a sparse family if 
the following properties hold:
\begin{enumerate}
\item
for each fixed $k\in\N_0$, 
the cubes $\{Q^k_j\}_{j\in J_k}$ 
are disjoint;
\item
if $\Omega_k:=\bigcup_{j\in J_k}Q^k_j$, 
then 
$\Omega_{k+1}\subset \Omega_k$;
\item
$|\Omega_{k+1}\cap Q^k_j|
\leq
\frac{1}{2}|Q^k_j|
$
for all $j\in J_k$.
\end{enumerate}
%Our aim in the next topic is to generalize Theorem \ref{th-151129-1} 
%and prove the local sharp maximal inequality on the weighted Morrey spaces
%as an application.

Now, we state our results related to the local sharp maximal operator. 
Throughout this paper, 
for $w\in A_\infty$,
we fix the parameter $\lambda_w'$ determined by $w$ 
so that 
$\lambda_w'<2^{-1-2^{n+3}[w]_{A_\infty}}$ 
and let 
$\lambda_w:=2^{-n-2}\lambda_w'$,
where
\[
[w]_{A_\infty}:=\lim_{r \uparrow \infty}[w]_{A_r}.
\]
\begin{theorem}\label{th-151129-3}
Let $0<s\leq q\leq p<\infty$ and $w\in A_\infty$.
\begin{enumerate}
\item
(Komori-Shirai type)
For all $f\in L^0(\R^n)$, it holds that
\begin{equation*}
\|f\|_{\mathcal{M}^p_q(w,w)}
\sim_{p,s,q,w}
\left\|M^{\sharp,d}_{\lambda_w}f\right\|_{\mathcal{M}^p_q(w,w)}
+
\|f\|_{\mathcal{M}^p_s(w,w)}.
\end{equation*}
\item
(Samko type)
Additionally, we assume $w\in\mathcal{B}_{p,q}$.
Then for all $f\in L^0(\R^n)$, it holds that
\begin{eqnarray}\label{160125-1}
%\lefteqn{
\|f\|_{\mathcal{M}^p_q(dx,w)}
%}
&\sim_{p,s,q,w}&
\left\|M^{\sharp,d}_{\lambda_w}f\right\|_{\mathcal{M}^p_q(dx,w)}\nonumber\\
&&+
\sup_{Q\in\mathcal{Q}}
\Phi_{p,q,w}(Q)
\left(\frac{1}{w(Q)}\int_Q|f(x)|^sdw(x)\right)^\frac{1}{s}.
\end{eqnarray}
\end{enumerate}
\end{theorem}

The last term
in the right-hand side of (\ref{160125-1})
is not so artificial
in view of \cite[(44)]{SaSh08}.
By the H\"{o}lder inequality,
it is easy to see
\begin{equation}\label{160121-10}
\sup_{Q\in\mathcal{Q}}
\Phi_{p,q,w}(Q)
\left(\frac{1}{w(Q)}\int_Q|f(x)|^sdw(x)\right)^\frac{1}{s}
\le 
\|f\|_{\mathcal{M}^p_q(dx,w)}
\end{equation}
for all $f \in L^0({\mathbb R}^n)$.
This corresponds to the relation:
\begin{equation}\label{160121-11}
\|f\|_{\mathcal{M}^p_s(w,w)}
\le
\|f\|_{\mathcal{M}^p_q(w,w)}
\end{equation}
for all $f \in L^0({\mathbb R}^n)$.

In the original sharp maximal inequality, 
we need to assume some integrability for $f$; 
see Theorem \ref{th-151201-5} for example.
However, in Theorem \ref{th-151129-3}, 
we can apply the sharp maximal inequality for any measurable functions.
Note that the same refinement for the unweighted case 
(but the inhomogeneous setting) 
is obtained in 
\cite[Theorem 1.3]{SaTa07-2}; 
see Theorem \ref{th-151201-4} in the present paper 
for the precise formulation.

If we impose a suitable condition for $f$, 
then we can recover the original version of the sharp maximal inequality, 
see Corollaries \ref{cr-151201-2} and \ref{cr-151204-1} as well.
\begin{theorem}\label{th-151129-4}
Let $0<s\leq q\leq p<\infty$, 
$w\in A_\infty$. 
\begin{enumerate}
\item
(Komori-Shirai type)
Assume that $f\in L^0(\R^n)$ satisfies
\begin{equation}\label{160121-2}
m_f(2^l Q)\to 0
\end{equation} 
as $l\to\infty$ for any $Q\in\mathcal{Q}$ 
and for some medians $\{m_f(2^lQ)\}_{l\in\N_0}$.
Then we have that
\begin{equation*}
\|f\|_{\mathcal{M}^p_s(w,w)}
\lesssim_{p,s,w}
\left\|M^{\sharp,d}_{\lambda_w}f\right\|_{\mathcal{M}^p_s(w,w)}
\leq
\left\|M^{\sharp,d}_{\lambda_w}f\right\|_{\mathcal{M}^p_q(w,w)}.
\end{equation*}
\item
(Samko type)
We assume
the weighted integral condition
(\ref{weight-integral}) for $p,q$ and $w$. 
Then for any $f\in L^0(\R^n)$ satisfying 
(\ref{160121-2}) for any $Q\in\mathcal{Q}$ 
and for some medians $\{m_f(2^lQ)\}_{l\in\N}$,
we have that
\begin{equation}\label{151215-1}
%\sup_{Q\in\mathcal{Q}}
%\Phi_{p,q,w}(Q)
%\left(\frac{1}{w(Q)}\int_Q|f(x)|^sdw(x)\right)^\frac{1}{s}
%\leq
\|f\|_{\mathcal{M}^p_q(dx,w)}
\lesssim_{p,q,w}
\left\|M^{\sharp,d}_{\lambda_w}f\right\|_{\mathcal{M}^p_q(dx,w)}.
\end{equation}
\end{enumerate}
\end{theorem}

\begin{remark}
In the case of the weighted Morrey space $\mathcal{M}^p_q(dx,w)$ 
of Samko type, 
the weighted integral condition
(\ref{weight-integral}) for $p,q,$ and $w$
and the decay condition
(\ref{160121-2})
are natural in view of the following example:
In fact, if we take 
$0<q<p<\infty$, 
$f_0\equiv1$ and 
$w_0(x):=|x|^{-\frac{q}{p}n}\in A_\infty\cap \mathcal{B}_{p,q}$,
then we know that 
$Mf_0\equiv1$, 
$f_0^\sharp=0$ 
and 
$0<\|Mf_0\|_{\mathcal{M}^p_q(dx,w_0)}
=
\||x|^{-\frac{n}{p}}\|_{\mathcal{M}^p_q(dx,dx)}
<\infty$,
which implies that 
the inequality 
\begin{equation*}
\|f\|_{\mathcal{M}^p_q(dx,w_0)}
\lesssim
\|f^\sharp\|_{\mathcal{M}^p_q(dx,w_0)}
\end{equation*}
fails even when $f$ satisfies
$Mf\in \mathcal{M}^p_q(dx,w_0)$.
One can notice that for $w_0(x)=|x|^{-\frac{q}{p}n}$, 
the weighted integral condition
(\ref{weight-integral}) for $p,q$ and $w_0$ fails to hold; 
see \cite{NS15}.  
\end{remark}

We note that for $0<s\leq q\leq p<\infty$,
\begin{equation*}
\sup_{Q\in\mathcal{Q}}
\Phi_{p,q,w}(Q)
\left(\frac{1}{w(Q)}\int_Q|f(x)|^sdw(x)\right)^\frac{1}{s}
\leq
\|f\|_{\mathcal{M}^p_q(dx,w)}
\end{equation*}
holds by H\"{o}lder's inequality.
Hence, combining Theorems \ref{th-151129-3} and \ref{th-151129-4}, 
we obtain the following:
\begin{corollary}\label{cr-151201-2}
Let $0< q\leq p<\infty$, 
$w\in A_\infty$. 
Assume that $f\in L^0(\R^n)$ satisfies
(\ref{160121-2})
for any $Q\in\mathcal{Q}$ 
and for some medians $\{m_f(2^lQ)\}_{l\in\N_0}$.
Then 
\begin{equation*}
\|f\|_{\mathcal{M}^p_q(w,w)}
\sim_{p,q,w}
\left\|M^{\sharp,d}_{\lambda_w}f\right\|_{\mathcal{M}^p_q(w,w)}.
\end{equation*}
Moreover, if
the weighted integral condition
(\ref{weight-integral}) for $p,q$ and $w$ holds, 
then we have that
\begin{equation*}
\|f\|_{\mathcal{M}^p_q(dx,w)}
\sim_{p,q,w}
\left\|M^{\sharp,d}_{\lambda_w}f\right\|_{\mathcal{M}^p_q(dx,w)}.
\end{equation*}
\end{corollary}

The remaining part of this paper is organized as follows:
In Section \ref{s2}, we collect some important estimates. 
%In Section \ref{s3}, we prove Theorems \ref{th-151215-1} and 
%\ref{th-151221-2} related to the boundedness 
%of the Hardy-Littlewood maximal operator.
%The key tool is the sparse family of cubes
%introduced in \cite{Hytonen11}.
Section \ref{s4} is devoted to the proof of Theorems \ref{th-151129-3} and
\ref{th-151129-4} related to the sharp maximal inequalities.
Here the decomposition by Lerner \cite{Lerner13-2} plays
an important role.
In Section \ref{s5}, we prove Theorems \ref{th-160626-1} and
\ref{th-160113-1} and then
give another application of the sharp maximal inequalities.
We take up the boundedness of the commutators.
%Section \ref{s6} considers some relations to the Morrey-Campanate spaces.
%In Section \ref{s7}, we give the proof of some fundamental lemmas 
%for the sake of completeness.

\section{Preliminaries}
\label{s2}

We first recall an important property of the class $\mathcal{B}_{p,q}$ and 
the relation between  \eqref{151221-1} and $M$. 
\begin{lemma}[\cite{NS15,Ta15}]\label{lm-160626-3}
Let $0<q\leq p<\infty$ and $w$ be a weight. 
\begin{enumerate}
\item
Then $w\in\mathcal{B}_{p,q}$ if and only if 
$\|\chi_{Q_0}\|_{\mathcal{M}^p_q(dx,w)}
\sim\Phi_{p,q,w}(Q_0)$ 
for any $Q_0\in\mathcal{Q}$;
see $(\ref{160627-1})$.
\item
We assume $q>1$ further. 
If $M$ is bounded on $\mathcal{M}^p_q(dx,w)$, 
then 
\eqref{151221-1} holds. 
\end{enumerate}
\end{lemma}

We also prove an important relation 
combining the $\mathcal{B}_{p,q}$ condition 
and (\ref{151221-1}) as follows:
\begin{lemma}\label{lm-160113-1}
Let $0<q\leq p<\infty$ and $w$ be a weight.
\begin{enumerate}
\item\label{aa}
$w$ satisfies (\ref{151221-1}) for all $Q\in\mathcal{Q}$, 
if and only if
$w\in\mathcal{B}_{p,q}$ and 
\begin{equation}\label{160113-4}
\frac{1}{|Q|}\Phi_{p,q,w}(Q)\left\|w^{-\frac{1}{q}}\chi_Q\right\|_{H^{q',n(1-q/p)}}
\leq
C_0
\end{equation}
holds for all $Q\in\mathcal{Q}$.
\item\label{bb}
For fixed $Q\in\mathcal{Q}$,
$w$ satisfies (\ref{160113-4}), 
if and only if 
\begin{equation}\label{160113-5}
\frac{1}{|Q|}\int_Qf(x)dx\leq
C_0
\frac{\|f\cdot\chi_Q\|_{\mathcal{M}^p_q(dx,w)}}
{\Phi_{p,q,w}(Q)}
\end{equation}  
holds for any $f\geq0$.
Here, the constant $C_0$ in (\ref{160113-5}) is the same as 
the one of (\ref{160113-4}).
Moreover, once $w$ satisfies 
(\ref{160113-5}) for all $Q$ with $C_0$ independent of $Q$, 
then $w$ is in the class $A_{q+1}$.
In particular, $w$ is a doubling weight.
\end{enumerate}
\end{lemma}

\begin{remark}\label{rm-160628-1}
From Lemma \ref{lm-160113-1}, particularly, we see that the inequality \eqref{151221-1}  
implies $w\in\mathcal{B}_{p,q}\cap A_{q+1}$. 
\end{remark}

\begin{proof}
To show the assertion \eqref{aa}, 
we have only to prove $w\in\mathcal{B}_{p,q}$ by 
imposing (\ref{151221-1}) for all $Q\in\mathcal{Q}$.
The converse is clear from the definition of $\Phi_{p,q,w}$.
To this end, we recall the notion of the sparse family.
Fix any $Q_0\in\mathcal{D}$ and let us 
show that 
\begin{equation*}
\|\chi_{Q_0}\|_{\mathcal{M}^p_q(dx,w)}
\sim_{n,p,q,w}
\Phi_{p,q,w}(Q_0).
\end{equation*}
Set $\gamma_0:=\frac{w(Q_0)}{|Q_0|}$ and 
take $a \gg 2^n$.
Define $\mathcal{D}_0:=\{Q_0\}$ 
and 
\begin{equation}\label{151215-2}
\mathcal{D}_k
:=
\left\{
Q\in\mathcal{D}(Q_0):
\frac{w(Q)}{|Q|}>a^k\gamma_0
\right\}
\quad
(k\in\N).
\end{equation}
By $\mathcal{D}_k^*:=\{Q^k_j\}_{j\in J_k}$, 
we denote the maximal subset of 
$\mathcal{D}_k$.
The maximality implies that 
for each $Q^k_j\in \mathcal{D}_k^*$, 
\begin{equation}\label{151212-1}
a^k\gamma_0
<\frac{w(Q^k_j)}{|Q^k_j|}
\leq
2^na^k\gamma_0.
\end{equation}
Then we claim that 
the family 
$\{Q^k_j\}_{k\in\N_0,i\in J_k}$ 
is a sparse family.
By the maximality, 
it is clear that 
$\{Q^k_j\}_{j\in J_k}$ 
is pairwise disjointed.
Moreover, 
for any $Q^{k+1}_j\in\mathcal{D}_{k+1}^*$, 
it follows from the definition of ${\mathcal D}_{k+1}$
that 
$\frac{w(Q^{k+1}_j)}{|Q^{k+1}_j|}>a^k\gamma_0$ by $a>1$, 
which implies $Q^{k+1}_j\in\mathcal{D}_k$.
Namely, we have 
$\Omega_{k+1}\subset 
\Omega_k$, 
where we denote 
\begin{equation}\label{160121-6}
\Omega_k:=\bigcup_{j\in J_k}Q^k_j.
\end{equation}
We also have that 
\begin{equation}\label{151215-3}
|Q^k_j\cap \Omega_{k+1}|
\leq
\frac{2^n}{a}|Q^k_j|,
\end{equation}
for any $k\in\N_0$ and $j\in J_k$.
Indeed, since $\{Q^{k+1}_i\}_{i\in J_{k+1}}$ is disjoint, 
thanks to
(\ref{151212-1})
\begin{eqnarray*}
|Q^k_j\cap \Omega_{k+1}|
&=&
\sum_{i\in J_{k+1}:Q^{k+1}_j\subset Q^k_j}
|Q^{k+1}_i|\\
&\leq&
\sum_{i\in J_{k+1}:Q^{k+1}_j\subset Q^k_j}
\frac{1}{a^{k+1}\gamma_0}w(Q^{k+1}_i)\\
&\leq&
\frac{1}{a^{k+1}\gamma_0}
w(Q^k_j)
\leq
\frac{2^n}{a}|Q^k_j|.
\end{eqnarray*}

%Since $\{Q^k_j\}_{k\in\N_0,j\in J_k}$ is a sparse family, 
%it follows that 
%\begin{equation*}
%\chi_{Q^k_j}(x)
%\leq
%C_{a,n}
%M[\chi_{Q^k_j\setminus \Omega_{k+1}}](x)
%\end{equation*}
%for any $k\in\N_0$ and $j\in J_k$.
Now, it follows from 
(\ref{151215-3}) that 
$
\frac{|Q^k_j\setminus \Omega_{k+1}|}{|Q^k_j|}
\geq
2^{-1}$,
for all $k\in\N_0$ and $j\in J_k$.
In particular, we focus on
the case of $k=0$ 
to obtain that 
$\frac{1}{|Q_0|}\int_{Q_0}\chi_{Q_0\setminus\Omega_1}(x)dx\geq2^{-1}$.
Using the dual inequality (\ref{dual}) and $a \gg 1$, we see that
\begin{equation*}
\frac{1}{2}
\leq
\frac{1}{|Q_0|}
\|\chi_{Q_0\setminus\Omega_1}\|_{\mathcal{M}^p_q(dx,w)}
\left\|w^{-\frac{1}{q}}\chi_{Q_0}\right\|_{H^{q',n(1-q/p)}}.
\end{equation*}
Moreover, 
it follows from  the assumption (\ref{151221-1}) that 
\begin{equation*}
\|\chi_{Q_0}\|_{\mathcal{M}^p_q(dx,w)}
\leq
C
\|\chi_{Q_0\setminus\Omega_1}\|_{\mathcal{M}^p_q(dx,w)}.
\end{equation*}
Now, let us recall that
\begin{eqnarray*}
&&\|\chi_{Q_0\setminus\Omega_1}\|_{\mathcal{M}^p_q(dx,w)}
=
\sup_{R\in\mathcal{D}(Q_0)}
|R|^\frac{1}{p}
\left(\frac{w(R\setminus\Omega_1)}{|R|}\right)^\frac{1}{q},\\
&&\Omega_1
=
\bigcup_{R\in\mathcal{D}_1}R,
\quad
\mathcal{D}_1
=
\left\{
R\in\mathcal{D}(Q_0):
\frac{w(R)}{|R|}>a\frac{w(Q_0)}{|Q_0|}
\right\}.
\end{eqnarray*}
This implies that 
for any $R\in\mathcal{D}(Q_0) \setminus \mathcal{D}_1$, 
$w(R\setminus\Omega_1)=0$ 
and hence, 
\begin{equation*}
\|\chi_{Q_0\setminus\Omega_1}\|_{\mathcal{M}^p_q(dx,w)}
=
\sup_{R\in\mathcal{D}(Q_0)\setminus\mathcal{D}_1}
|R|^\frac{1}{p}
\left(\frac{w(R\setminus\Omega_1)}{|R|}\right)^\frac{1}{q}
\leq
a^\frac{1}{q}\Phi_{p,q,w}(Q_0).
\end{equation*}
In summary, 
we obtain that
\begin{equation*}
\|\chi_{Q_0}\|_{\mathcal{M}^p_q(dx,w)}
\leq
C_{a,n,p,q,w}
\Phi_{p,q,w}(Q_0).
\end{equation*}
Meanwhile, since the converse inequality is trivial, 
it follows $\|\chi_{Q_0}\|_{\mathcal{M}^p_q(dx,w)}\sim\Phi_{p,q,w}(Q_0)$ 
and hence $w\in\mathcal{B}_{p,q}$.

Next, let us show the assertion \eqref{bb}.
Once we suppose (\ref{160113-4}), then 
we deduce (\ref{160113-5}) immediately
from  (\ref{dual}).
Conversely, if we assume (\ref{160113-5}), 
then 
we employ (\ref{dual2}) to obtain 
\begin{eqnarray*}
\left\|w^{-\frac{1}{q}}\chi_Q\right\|_{H^{q',n(1-q/p)}}
&=&
\sup_{f\geq0:\|f\chi_Q\|_{\mathcal{M}^p_q(dx,w)}\leq1}
\int_Qf(x)dx\\
&\leq&
\sup_{f\geq0:\|f\chi_Q\|_{\mathcal{M}^p_q(dx,w)}\leq1}
C_0\frac{|Q|}{\Phi_{p,q,w}(Q)}\|f\chi_Q\|_{\mathcal{M}^p_q(dx,w)}\\
&=&
C_0\frac{|Q|}{\Phi_{p,q,w}(Q)},
\end{eqnarray*} 
which implies (\ref{160113-4}).

Finally, we shall show that 
$w$ is in $A_{q+1}$
%satisfies the doubling condition 
by imposing 
(\ref{160113-5}) for all $Q\in\mathcal{Q}$.
Note that 
\begin{equation*}
\left\|w^{-\frac{1}{q}}\chi_Q\right\|_{\mathcal{M}^p_q(dx,w)}
=
\|\chi_Q\|_{\mathcal{M}^p_q(dx,dx)}
=
|Q|^\frac{1}{p}.
\end{equation*}
With this in mind, 
by putting $f:=w^{-\frac{1}{q}}$ in (\ref{160113-5}), 
then it follows that 
\begin{equation*}
\frac{1}{|Q|}\int_Qw(x)^{-\frac{1}{q}}dx
\leq
C_0
\frac{|Q|^\frac{1}{p}}{\Phi_{p,q,w}(Q)}
=
C_0
\left(\frac{|Q|}{w(Q)}\right)^\frac{1}{q},
\end{equation*}
which implies that 
$[w]_{A_{q+1}}\leq C_0$.
%Another proof for that (\ref{160113-5}) implies $w$ is doubling%%%%%%%%%%%%%%%%%%%%%%%%%%%%
\if0
Our aim is to show that $w(2Q_0)\leq Cw(Q_0)$ holds 
for any $Q_0\in\mathcal{Q}$.
To this end, we fix any $Q_0\in\mathcal{Q}$ and 
consider the Calder\'{o}n-Zygmund decomposition as in the proof of 
the assertion \ref{aa} to
construct $\Omega_1\subset Q_0$ such that 
\begin{equation}\label{160114-6}
\frac{1}{|Q_0|}|Q_0\setminus\Omega_1|
\geq
\frac{1}{2},
\quad
\|\chi_{Q_0\setminus\Omega_1}\|_{\mathcal{M}^p_q(dx,w)}
\leq
2^{\frac{1}{q}(n+1)}\Phi_{p,q,w}(Q_0).
\end{equation}
Now, we take 
$Q=2Q_0$ and $f=\chi_{Q_0\setminus\Omega_1}$ in (\ref{160113-5}). 
Then it follows that 
\begin{equation*}
\frac{1}{|2Q_0|}|Q_0\setminus\Omega_1|
\leq
C_0
\frac{\|\chi_{Q_0\setminus\Omega_1}\|_{\mathcal{M}^p_q(dx,w)}}
{\Phi_{p,q,w}(2Q_0)}.
\end{equation*}
Using (\ref{160114-6}), 
we see that 
$\Phi_{p,q,w}(2Q_0)
\leq
2^{(1/q+1)(n+1)}\Phi_{p,q,w}(Q_0)$, 
or equivalently 
$w(2Q_0)
\leq
2^{(1/q+1)(n+1)}
2^{n\left(\frac{1}{q}-\frac{1}{p}\right)}
C_0
w(Q_0).
$
By putting
$C:=2^{(1/q+1)(n+1)}
2^{n\left(\frac{1}{q}-\frac{1}{p}\right)}
C_0$, 
we obtain the conclusion.
\fi
%%%%%%%%%%%%%%%%%%%%%%%%%%%%%%%%%%%%%%%%%%%%%%%%%%%%%%%%%%%%
\end{proof}

In \cite{Nakai94}, Nakai quantitatively showed the self-improvement 
of the integral condition; see \cite[Proposition 2.6]{NNS15} as well.
\begin{proposition}[\cite{Nakai94}]\label{pr-160115-1}
Assume that a positive function $f:(0,\infty)\to(0,\infty)$
satisfies the integral condition:
\begin{equation}\label{160115-1}
\int_1^\infty \frac{1}{f(rs)}\frac{ds}{s}
\leq
\frac{C_0}{f(r)}
\quad
(r>0)
\end{equation}
for some constant $C_0>0$.
Then for all $\delta\in(0,C_0^{-1})$, 
we have 
\begin{equation}\label{160115-2}
\int_1^\infty \frac{s^\delta}{f(rs)}\frac{ds}{s}
\leq
\frac{C_0}{1-C_0\delta}\frac{1}{f(r)}
\quad
(r>0).
\end{equation}
\end{proposition}
Using Proposition \ref{pr-160115-1}, 
we can show the self-improvement of the
weighted integral condition.
We will invoke the following observation in the
proof of Theorem \ref{th-160113-1}.
\begin{lemma}\label{lm-160115-1}
Let $0<q\leq p<\infty$ and $w$ satisfies the 
weighted integral condition (\ref{weight-integral}) for 
$p,q$ and $w$.
Then we have 
\begin{equation}\label{160115-3}
\int_1^\infty
\frac{{\rm log}s}{\Phi_{p,q,w}(sQ)}
\frac{ds}{s}
\leq
\frac{C}{\Phi_{p,q,w}(Q)}
\quad
(Q\in\mathcal{Q}),
\end{equation}
or equivalently, 
\begin{equation}\label{160115-4}
\sum_{k=1}^\infty\frac{k}{\Phi_{p,q,w}(2^kQ)}
\leq
\frac{C}{\Phi_{p,q,w}(Q)}
\quad
(Q\in\mathcal{Q}).
\end{equation}
\end{lemma}

\begin{proof}
Fix any cube $Q$.
We put 
$f(s):=\Phi_{p,q,w}(sQ)$ for $s\geq1$ and
$f(s):=\Phi_{p,q,w}(Q)$ for $s\leq1$.
Since we assume the weighted integral condition 
(\ref{weight-integral})
for $p,q$ and $w$, 
the integral condition
for the above $f$ holds.
Hence, we see that 
by Proposition \ref{pr-160115-1}, 
\begin{equation*}
\int_1^\infty
\frac{s^\delta}{\Phi_{p,q,w}(sQ)}
\frac{ds}{s}
\leq
\frac{C}{\Phi_{p,q,w}(Q)}
\end{equation*}
for some $\delta>0$.
If we notice that 
${\rm log}s\leq C_\delta s^\delta$, 
then we obtain (\ref{160115-3}).
Note that the equivalence between 
(\ref{160115-3}) and (\ref{160115-4}) 
is clear from this observation.
\end{proof}

Next, we collect the fundamental properties of the median and oscillation;
see \cite{Hytonen11} for the details.
Let $f\in L^0(\R^n)$, $\lambda\in(0,2^{-1})$ and $Q\in\mathcal{Q}$.
For any median $m_f(Q)$, we have that
\begin{equation}\label{151130-1}
|m_f(Q)|
\leq
\left(f\cdot\chi_Q\right)^*(\lambda|Q|),
\quad
\lim_{\ell(Q)\to0,Q\ni x}m_f(Q)
=
f(x)
\quad
({\rm a.e.}\ x\in\R^n).
\end{equation}
See \cite[Lemma 2.2]{Fujii91} for the details.
Moreover, for the oscillation of $f$, we have that
\begin{equation}\label{151130-2}
\omega_\lambda(f;Q)
\leq
\left((f-m_f(Q))\chi_Q\right)^*(\lambda|Q|)
\leq
2\omega_\lambda(f;Q).
\end{equation}

Now, we introduce the notion of $w$-sparse family for $w\in A_\infty$.
Let $w\in A_\infty$ and $\lambda_w'<2^{-1-2^{n+3}[w]_{A_\infty}}$.
We say that 
the family $\{Q^k_j\}_{k\in\N_0,j\in J_k}$ 
is a $w$-sparse family if 
the following properties hold:
\begin{enumerate}
\item
for each fixed $k\in\N_0$, 
the cubes $\{Q^k_j\}_{j\in J_k}$ 
are disjoint;
\item
if $\Omega_k:=\bigcup_{j\in J_k}Q^k_j$, 
then 
$\Omega_{k+1}\subset \Omega_k$;
\item
$|\Omega_{k+1}\cap Q^k_j|
\leq
\lambda_w'|Q^k_j|
$
for all $j\in J_k$.
\end{enumerate}
The only difference from the original definition of the sparse family 
is the constant $\lambda_w'$ appearing in the condition $(3)$.
That is, if we replace the constant 
$\lambda_w'$ by $2^{-1}$ in the condition $(3)$, 
then the 
$w$-sparse family turns into the original sparse family introduced in \cite{Hytonen11}.
The advantage of $w$-sparse family is as follows:
\begin{lemma}\label{lm-151201}
Let $w\in A_\infty$, $\lambda_w'<2^{-1-2^{n+3}[w]_{A_\infty}}$ 
and 
$\{Q^k_j\}_{k\in\N_0,j\in J_k}$ be a $w$-sparse family.
Then there exists a constant $C_w>1$ such that 
\begin{equation*}
w(Q^k_j)
\leq
C_w
w(Q^k_j\cap \Omega_{k+1}^c),
\end{equation*}
holds for all 
$k\in\N_0$ and $j\in J_k$.
Here, the constant $C_w$ can be written by 
\begin{equation*}
C_w=\left(1-2(\lambda_w')^{\frac{\varepsilon}{1+\varepsilon}}\right)^{-1},
\quad
\varepsilon:=
\frac{1}{2^{n+3}[w]_{A_\infty}}.
\end{equation*}
\end{lemma}

We can find such a property in many papers; see for example \cite{GR85}, 
and hence, we omit the proof. 
We will invoke Lemma \ref{lm-151201} to show 
Theorem \ref{th-151129-2} which is a crucial part of Theorem \ref{th-151129-3}.

%%%%%%%%%%%%%%%%%%%%%%%%%%%%%%%%%%%%%
\if0
\begin{proof}
We can find such property in many papers; for example \cite{GR85}, 
but we give the complete proof.
The proof is based on the reverse H\"{o}lder inequality.
Since $\{Q^k_j\}_{k\in\N_0,j\in J_k}$ is a $w$-sparse family,
we have that 
\begin{eqnarray*}
w(Q^k_j\cap\Omega_{k+1})
&\leq&
\left(\int_{Q^k_j}w(x)^{1+\varepsilon}dx\right)^\frac{1}{1+\varepsilon}
\cdot
\left|\Omega_{k+1}\cap Q^k_j\right|^{\frac{\varepsilon}{1+\varepsilon}}\\
&\leq&
(\lambda_w')^\frac{\varepsilon}{1+\varepsilon}
\left(\int_{Q^k_j}w(x)^{1+\varepsilon}dx\right)^\frac{1}{1+\varepsilon}
\cdot
\left|Q^k_j\right|^{\frac{\varepsilon}{1+\varepsilon}}.
\end{eqnarray*}
Hence, by the reverse H\"{o}lder inequality, 
we see that 
$w(Q^k_j\cap \Omega_{k+1})
\leq
2(\lambda_w')^\frac{\varepsilon}{1+\varepsilon}w(Q^k_j)$.
Here, note that 
the assumption: 
$\lambda_w'<2^{-1-2^{n+3}[w]_{A_\infty}}$
ensures  
$2(\lambda_w')^\frac{\varepsilon}{1+\varepsilon}<1$.
With the trivial equality: 
$w(Q^k_j)=w(Q^k_j\cap\Omega_{k+1}^c)
+
w(Q^k_j\cap\Omega_{k+1})$ in mind, 
it follows that 
\begin{equation*}
w(Q^k_j)
\leq
\left(1-2(\lambda_w')^\frac{\varepsilon}{1+\varepsilon}\right)^{-1}
w(Q^k_j\cap\Omega_{k+1}^c).
\end{equation*}
\end{proof}
\fi
%%%%%%%%%%%%%%%%%%%%%%%%%%%%%%%%%%%

Here, we give an observation combining our results and 
the original sharp maximal inequalities.
We will employ the following lemma to show Corollary \ref{cr-151204-1}. 
\begin{lemma}\label{lm-151201-3}
Let $f\in L^0(\R^n)$.
Assume that $Mf\in\mathcal{M}^{p_0}_{q_0}(v_0,w_0)$ 
for some $0<q_0\leq p_0<\infty$ and some weights $v_0,w_0$ 
satisfying that 
\begin{equation}\label{151203-1}
\lim_{l\to\infty}\Phi_{p_0,q_0,v_0,w_0}(2^lQ)
=\infty,
\quad
(Q\in\mathcal{Q}),
\quad
\Phi_{p_0,q_0,v_0,w_0}(Q)
:=
v_0(Q)^{\frac{1}{p_0}-\frac{1}{q_0}}w_0(Q)^\frac{1}{q_0}.
\end{equation}
Then for any $Q\in\mathcal{Q}$ and any medians $\{m_f(2^lQ)\}_{l\in\N_0}$, 
it holds that 
\begin{equation}\label{151203-2}
\lim_{l\to\infty}m_f(2^lQ)=0.
\end{equation}
In particular, we have the following:
\begin{enumerate}
\item
If we assume that $Mf\in \mathcal{M}^{p_0}_{q_0}(w_0,w_0)$ 
for some $0<q_0\leq p_0<\infty$ and $w_0\in A_\infty$.
Then for any $Q\in\mathcal{Q}$ and any medians $\{m_f(2^lQ)\}_{l\in\N_0}$, 
(\ref{151203-2}) holds.
\item
If we assume that $Mf\in \mathcal{M}^{p_0}_{q_0}(dx,w_0)$ 
for some $0<q_0\leq p_0<\infty$ and $w_0\in A_\infty$ satisfying 
the weighted integral condition
(\ref{weight-integral}) for $p_0,q_0$ and $w_0$.
Then for any $Q\in\mathcal{Q}$ and any medians $\{m_f(2^lQ)\}_{l\in\N_0}$, 
(\ref{151203-2}) holds.
\end{enumerate}
\end{lemma}

\begin{proof}
In view of (\ref{151130-1}), 
for $\lambda\in(0,2^{-1})$, 
we have that 
\begin{equation*}
|m_f(2^lQ)|
%detail%%%%%%%%%%%%%%%%%%%%%%%%%%%%%%%%%%%
\leq
\left(f\cdot\chi_{2^lQ}\right)^*(\lambda|2^lQ|)
%&=&
%\frac{1}{\lambda|2^lQ|}\int_0^{\lambda|2^lQ|}
%\left(f\cdot\chi_{2^lQ}\right)^*(\lambda|2^lQ|)dt\\
\leq
\frac{1}{\lambda|2^lQ|}\int_0^\infty
\left(f\cdot\chi_{2^lQ}\right)^*(t)dt
%&=&
%\frac{1}{\lambda|2^lQ|}
%\int_{2^lQ}|f(y)|dy\\
\leq
\frac{1}{\lambda}
\inf_{x\in2^lQ}Mf(x).
\end{equation*}
With the definition of the norm in mind, 
by taking $\lambda=4^{-1}$ for example, 
we notice that 
\begin{equation}\label{151203-3}
|m_f(2^lQ)|
\lesssim
\frac{1}{w_0(2^lQ)^\frac{1}{q_0}}
\left(\int_{2^lQ}Mf(x)^{q_0}dw_0(x)\right)^\frac{1}{q_0}
\leq
\frac{
\|Mf\|_{\mathcal{M}^{p_0}_{q_0}(v_0,w_0)}
}{\Phi_{p_0,q_0,v_0,w_0}(2^lQ)}.
\end{equation}
By our assumption, 
$\|Mf\|_{\mathcal{M}^{p_0}_{q_0}(v_0,w_0)}<\infty$.
Hence, the assumption (\ref{151203-1}) yields the conclusion (\ref{151203-2}).

Particularly, let us consider the case of $v_0=w_0\in A_\infty$. 
The condition (\ref{151203-1}) reads
$
\lim_{l\to\infty}
w_0(2^lQ)^\frac{1}{p_0}
=\infty
$
for all $Q\in\mathcal{Q}$.
Note that 
$w_0\in A_\infty$ implies $w_0\notin L^1(dx)$;
see \cite[Chapter 1, \S 8.6(a)]{St2}.
In fact, by the reverse H\"{o}lder inequality, 
for any $R>0$, we have that
\begin{equation*}
\left(
\frac{1}{|Q(R)|}\int_{Q(R)}w_0(x)^{1+\varepsilon}dx
\right)^\frac{1}{1+\varepsilon}
\leq
\frac{2w_0(Q(R))}{|Q(R)|},
\end{equation*}
which implies 
\begin{equation*}
|Q(R)|^{1-\frac{1}{1+\varepsilon}}
\left(
\int_{Q(R)}w_0(x)^{1+\varepsilon}dx
\right)^\frac{1}{1+\varepsilon}
\leq
2w_0(Q(R)).
\end{equation*}
By tending $R\to\infty$, 
we obtain $\|w_0\|_{L^1(dx)}=\infty$.
Thus, we see that 
$\lim_{l\to\infty}
w_0(2^lQ)^\frac{1}{p}
=\|w_0\|_{L^1(dx)}^\frac{1}{p}
=\infty.
$

In the case of $v_0=dx$, 
the condition (\ref{151203-1}) reads  
$
\lim_{l\to\infty}
\Phi_{p_0,q_0,w_0}(2^lQ)
=\infty
$
for all $Q\in\mathcal{Q}$.
However, this condition follows directly from 
the weighted integral condition
(\ref{weight-integral}) for $p_0,q_0$ and 
$w_0$.
\end{proof}

\begin{remark}
The assumption (\ref{151203-1}) is essential in the following sense:
If the assumption (\ref{151203-1}) does not hold, 
then $Mf=f\equiv1$ is in $\mathcal{M}^{p_0}_{q_0}(v_0,w_0)$, 
but (\ref{151203-2}) fails for any $Q\in\mathcal{Q}$.
\end{remark}

To ensure the well definedness of 
the singular integral operators 
and commutators on weighted Morrey spaces 
in Theorems \ref{th-160626-1} and \ref{th-160113-1}, 
we need the self-improvement property of the boundedness of $M$ 
on $\mathcal{M}^p_q(dx,w)$.
The self-improvement property of $M$ on the Banach function spaces 
is obtained in \cite[Theorem 1.2]{Lerner-Ombrosi10}.
\begin{lemma}\label{lm-160113-2}
Let $1< q\leq p<\infty$ and $w$ be a weight.
Assume that 
the Hardy-Littlewood maximal operator $M$ is bounded on 
$\mathcal{M}^p_q(dx,w)$. 
Then there exists $r=r(p,q,w)>1$ such that 
$M^{(r)}$ is also bounded on $\mathcal{M}^p_q(dx,w)$.
\end{lemma}
Since the proof of Lemma \ref{lm-160113-2} is almost the same as 
the one of \cite[Theorem 1.2]{Lerner-Ombrosi10}, 
we omit the proof here. 
We can rephrase the assertion of 
Lemma \ref{lm-160113-2} as follows: 
Once we obtain the boundedness of $M$ on 
$\mathcal{M}^p_q(dx,w)$ for some 
$1< q\leq p<\infty$ and some weight $w$, 
then there exists a small $\varepsilon\in(0,1)$ 
such that 
$M$ is also bounded on 
$\mathcal{M}^{p(1-\varepsilon)}_{q(1-\varepsilon)}(dx,w)$.
This property yields the following corollary:
\begin{corollary}\label{cr-160113-4}
Assume that 
the Hardy-Littlewood maximal operator $M$ is bounded on 
$\mathcal{M}^p_q(dx,w)$ for some $1< q\leq p<\infty$ and 
some weight $w$.
Then for any $Q\in\mathcal{Q}$ and any $f\in \mathcal{M}^p_q(dx,w)$, 
$f\cdot\chi_Q$ is in $L^{\frac{1}{1-\varepsilon}}(\R^n)$. 
\end{corollary}
\begin{proof}
Let $\varepsilon$ be as above.
In this proof, we abbreviate $p(1-\varepsilon)$ and $q(1-\varepsilon)$ 
by $p_\varepsilon$ and $q_\varepsilon$ respectively.
By the above observation, 
%we can find small $\varepsilon\in(0,1)$ such that 
we know that
$M$ is bounded on 
$\mathcal{M}^{p_\varepsilon}_{q_\varepsilon}(dx,w)$.
By Lemmas \ref{lm-160626-3} and \ref{lm-160113-1}, 
it follows that 
\begin{equation*}
\frac{1}{|Q|}\Phi_{p,q,w}(Q)^\frac{1}{1-\varepsilon}
\left\|w^{-\frac{1}{q_\varepsilon}}\chi_Q\right\|_{H^{q_\varepsilon',n(1-q/p)}}
\lesssim
1,
\end{equation*}
for all $Q\in\mathcal{Q}$.
Here, we used a simple fact that 
$\Phi_{p_\varepsilon,q_\varepsilon,w}(Q)
=
\Phi_{p,q,w}(Q)^\frac{1}{1-\varepsilon}$.
With this in mind, 
by invoking the dual inequality (\ref{dual}), 
we have that
\begin{eqnarray*}
\int_{\R^n}|f(x)\chi_Q(x)|^\frac{1}{1-\varepsilon}dx
&\leq&
C
\left\||f|^\frac{1}{1-\varepsilon}w^{\frac{1}{q_\varepsilon}}\right\|_{\mathcal{M}^{p_\varepsilon}_{q_\varepsilon}(dx,dx)}
\left\|w^{-\frac{1}{q_\varepsilon}}\chi_Q\right\|_{H^{q_\varepsilon',n(1-q/p)}}\\
&\leq&
C
\left\||f|^\frac{1}{1-\varepsilon}\right\|_{\mathcal{M}^{p_\varepsilon}_{q_\varepsilon}(dx,w)}
\frac{|Q|}{\Phi_{p,q,w}(Q)^\frac{1}{1-\varepsilon}}.
\end{eqnarray*}
If we notice that 
$
\left\||f|^\frac{1}{1-\varepsilon}\right\|_{\mathcal{M}^{p_\varepsilon}_{q_\varepsilon}(dx,w)}
=
\|f\|_{\mathcal{M}^p_q(dx,w)}^\frac{1}{1-\varepsilon},
$
then we see that 
\begin{equation*}
\|f\cdot\chi_Q\|_{L^{\frac{1}{1-\varepsilon}}(\R^n)}
\leq
C
\frac{|Q|^{1-\varepsilon}}{\Phi_{p,q,w}(Q)}
\|f\|_{\mathcal{M}^p_q(dx,w)}<\infty.
\end{equation*}
\end{proof}

\section{Proof of Theorems \ref{th-151129-3} and \ref{th-151129-4}}
\label{s4}

Our first observation is to generalize Theorem \ref{th-151129-1} 
to the $w$-sparse family setting.
\begin{proposition}\label{pr-151201-1}
Let $w\in A_\infty$, $\lambda_w'<2^{-1-2^{n+3}[w]_{A_\infty}}$ 
and
$\lambda_w:=2^{-n-2}\lambda_w'$.
For $Q_0\in\mathcal{Q}$ and $f:Q_0\to\R$, 
there exists a $w$-sparse family 
$\{Q^k_j\}_{k\in\N_0,j\in J_k}\subset \mathcal{D}(Q_0)$
such that for a.e. $x\in\R^n$,
 \begin{equation*}
|f(x)-m_f(Q_0)|
\leq
4M^{\sharp,d}_{\lambda_w;Q_0}f(x)
+
2\sum_{k\in\N_0}\sum_{j\in J_k}\omega_{\lambda_w}(f;Q^k_j)\chi_{Q^k_j}(x).
\end{equation*}
\end{proposition}
Since we need only a slight modification of the proof of 
Theorem \ref{th-151129-1} to show Proposition \ref{pr-151201-1}, 
we omit the proof here. 
By employing Proposition \ref{pr-151201-1}, 
we can show the following local estimate:

\begin{theorem}{\label{th-151129-2}}
Let $0<q,s<\infty$, $w,v\in A_\infty$ 
and $Q_0\in\mathcal{Q}$.
We choose $\lambda_w'>0$ so that 
$\lambda_w'<2^{-1-2^{n+3}[w]_{A_\infty}}$
and let $\lambda_w:=2^{-n-2}\lambda_w'$.
Then for all $f\in L^0(\R^n)$, we have 
\begin{eqnarray*}
\left(\int_{Q_0}|f(x)|^qdw(x)\right)^\frac{1}{q}
&\lesssim_{q,s,w,v}&
\left(
\int_{Q_0}
M^{\sharp,d}_{\lambda_w;Q_0}f(x)^q
dw(x)
\right)^\frac{1}{q}\\
&&+
w(Q_0)^\frac{1}{q}
\left(\frac{1}{v(Q_0)}
\int_{Q_0}
|f(x)|^sdv(x)\right)^\frac{1}{s}.
\end{eqnarray*}
\end{theorem}

\begin{proof}
We take a median $m_f(Q_0)$ and use the quasi-triangle inequality 
to get
\begin{equation*}
\left(
\int_{Q_0}|f(x)|^qdw(x)
\right)^\frac{1}{q}
\lesssim_q
\left(
\int_{Q_0}|f(x)-m_f(Q_0)|^qdw(x)
\right)^\frac{1}{q}
+
w(Q_0)^\frac{1}{q}|m_f(Q_0)|.
\end{equation*}
For the first term, by applying Proposition \ref{pr-151201-1}, 
we obtain that
\begin{eqnarray*}
&&\left(
\int_{Q_0}|f(x)-m_f(Q_0)|^qdw(x)
\right)^\frac{1}{q}\\
&\lesssim_q&
\left(
\int_{Q_0}M^{\sharp,d}_{\lambda_w;Q_0}f(x)^qdw(x)
\right)^\frac{1}{q}
+
\left[\int_{Q_0}\left(\sum_{k\in\N_0}\sum_{j\in J_k}\omega_{\lambda_w}(f;Q^k_j)\chi_{Q^k_j}(x)\right)^q
dw(x)
\right]^\frac{1}{q}.
\end{eqnarray*}
Now, we focus on 
\begin{equation*}
{\rm I}
:=
\left[\int_{Q_0}\left(\sum_{k\in\N_0}\sum_{j\in J_k}\omega_{\lambda_w}(f;Q^k_j)\chi_{Q^k_j}(x)\right)^q
dw(x)
\right]^\frac{1}{q}.
\end{equation*}
We recall that $\{Q^k_j\}_{k\in\N_0,j\in J_k}$ is a $w$-sparse family, 
in particular, $w(Q^k_j)\leq C_ww(Q^k_j\cap\Omega_{k+1}^c)$ holds 
by Lemma \ref{lm-151201},
where $\Omega_{k+1}$ is given by
(\ref{160121-6}).
This implies that 
$\chi_{Q^k_j}(x)
\leq
C_w
M_w[\chi_{Q^k_j\cap\Omega_{k+1}^c}](x)$, 
where the weighted Hardy-Littlewood maximal operator $M_w$ is defined by 
\begin{equation*}
M_wf(x):=\sup_{Q\in\mathcal{Q}}\frac{\chi_Q(x)}{w(Q)}\int_Q|f(y)|dw(y).
\end{equation*}
Hence, by taking 
$\eta>\max{(1,q^{-1})}$
and employing the boundedness
of $M_w$ on $L^{\eta q}(\ell^{\eta})(w)$,
it follows that 
\begin{eqnarray*}
{\rm I}
&\lesssim_w&
\left[\int_{Q_0}\left(\sum_{k\in\N_0}\sum_{j\in J_k}\omega_{\lambda_w}(f;Q^k_j)
M_w[\chi_{Q^k_j\cap\Omega_{k+1}^c}](x)^\eta\right)^q
dw(x)
\right]^\frac{1}{q}\\
&\lesssim_{q,w}&
\left[\int_{\R^n}\left(\sum_{k\in\N_0}\sum_{j\in J_k}\omega_{\lambda_w}(f;Q^k_j)
\chi_{Q^k_j\cap\Omega_{k+1}^c}(x)\right)^q
dw(x)
\right]^\frac{1}{q}.
\end{eqnarray*} 
If we notice that 
$\omega_{\lambda_w}(f;Q^k_j)\chi_{Q^k_j}(x)
\leq
M^{\sharp,d}_{\lambda_w;Q_0}f(x)$ 
and that 
$\sum_{k\in\N_0}\sum_{j\in J_k}\chi_{Q^{k_j}\cap\Omega_{k+1}^c}\leq
\chi_{Q_0}$ by the disjointness of $\{Q^k_j\cap\Omega_{k+1}^c\}_{k\in\N_0,j\in J_k}$,
then we see that 
\begin{equation*}
{\rm I}
\lesssim_{q,w}
\left(\int_{Q_0}M^{\sharp,d}_{\lambda_w;Q_0}f(x)^qdw(x)
\right)^\frac{1}{q}.
\end{equation*}

Next, we evaluate the second term: 
$w(Q_0)^\frac{1}{q}|m_f(Q_0)|$.
Since $v\in A_\infty$, we can find 
$u\in(1,\infty)$ such that 
$v\in A_u$. 
We set $r:=\frac{s}{u}$.
By (\ref{151130-1}), for $\lambda\in(0,2^{-1})$, 
we have that 
\begin{eqnarray*}
|m_f(Q_0)|
&\leq&
\left(f\cdot\chi_{Q_0}\right)^*(\lambda|Q_0|)\\
&\leq&
\left(\frac{1}{\lambda|Q_0|}\int_0^{\lambda|Q_0|}
\left(f\cdot\chi_{Q_0}\right)^*(t)^rdt
\right)^\frac{1}{r}\\
&\lesssim_r&
\left(\frac{1}{|Q_0|}\int_{Q_0}|f(x)|^rdx
\right)^\frac{1}{r}.
\end{eqnarray*} 
Moreover, by $v\in A_u$ and $ru=s$, 
we see that 
\begin{eqnarray*}
|m_f(Q_0)|
&\lesssim_{r}&
[v]_{A_u}^\frac{1}{ru}
\left(\frac{1}{|Q_0|}\int_{Q_0}|f(x)|^{ru}dv(x)
\right)^\frac{1}{ru}
\cdot
\left(\frac{|Q_0|}{v(Q_0)}
\right)^\frac{1}{ru}\\
&=&
[v]_{A_u}^\frac{1}{s}
\left(\frac{1}{v(Q_0)}
\int_{Q_0}|f(x)|^sdv(x)
\right)^\frac{1}{s}.
\end{eqnarray*}
Thus, we complete the proof.
\end{proof}

As a corollary, we obtain the following:
\begin{corollary}\label{cr-151201-1}
Let $0<q\leq p<\infty$, $s>0$, $w_1$ be a weight, 
and $w_2,v\in A_\infty$.
Then for any $f\in L^0(\R^n)$, one has
\begin{eqnarray*}
\|f\|_{\mathcal{M}^p_q(w_1,w_2)}
&\lesssim_{q,s,w_2,v}&
\left\|M^{\sharp,d}_{\lambda_w}f\right\|_{\mathcal{M}^p_q(w_1,w_2)}\\
&&
+
\sup_{Q\in\mathcal{Q}}
w_1(Q)^\frac{1}{p}
\left(\frac{w_2(Q)}{w_1(Q)}\right)^\frac{1}{q}
\left(
\frac{1}{v(Q)}
\int_Q|f(x)|^sv(x)dx\right)^\frac{1}{s}.
\end{eqnarray*}
In particular, the following two types of weighted sharp maximal inequalities hold:
\begin{enumerate}
\item
(Komori-Shirai type: $w_1=w_2=v=w\in A_\infty$)
\begin{equation*}
\|f\|_{\mathcal{M}^p_q(w,w)}
\lesssim_{s,q,w}
\left\|M^{\sharp,d}_{\lambda_w}f\right\|_{\mathcal{M}^p_q(w,w)}
+
\|f\|_{\mathcal{M}^p_s(w,w)}.
\end{equation*}
\item
(Samko type: $w_1=dx$, $w_2=v=w\in A_\infty$)
\begin{eqnarray*}
\lefteqn{
\|f\|_{\mathcal{M}^p_q(dx,w)}
}\\
&\lesssim_{s,q,w}&
\left\|M^{\sharp,d}_{\lambda_w}f\right\|_{\mathcal{M}^p_q(dx,w)}
+
\sup_{Q\in\mathcal{Q}}
\Phi_{p,q,w}(Q)
\left(\frac{1}{w(Q)}\int_Q|f(x)|^sdw(x)\right)^\frac{1}{s}.
\end{eqnarray*}
\end{enumerate}
\end{corollary}

Next, we consider the converse inequality. 
To this end, 
we recall the sufficient conditions
for the boundedness of $M$
on weighted Morrey spaces.
\begin{theorem}\label{th-151202-2}
Let $1<q\leq p<\infty$ and $w\in A_q$.
\begin{enumerate}
\item
\cite{KoSh09}
We have that 
$
\|Mf\|_{\mathcal{M}^p_q(w,w)}
\lesssim_{p,q,w}
\|f\|_{\mathcal{M}^p_q(w,w)}
$
for all $f\in \mathcal{M}^p_q(w,w))$.
\item
\cite{NS15}
In addition, we assume $w\in \mathcal{B}_{p,q}$.
Then we have that
$
\|Mf\|_{\mathcal{M}^p_q(dx,w)}
\lesssim_{p,q,w}
\|f\|_{\mathcal{M}^p_q(dx,w)}
$
for all 
$f\in \mathcal{M}^p_q(dx,w)$.
\end{enumerate}
\end{theorem}

\begin{proposition}\label{pr-151201-2}
Let $0<q\leq p<\infty$, $\lambda\in(0,2^{-1})$ 
and $w\in A_\infty$.
Then for $f\in L^0(\R^n)$, 
we have that
\begin{equation}\label{160121-8}
\left\|M^{\sharp,d}_{\lambda}f\right\|_{\mathcal{M}^p_q(w,w)}
\lesssim_{\lambda,p,q,w}
\|f\|_{\mathcal{M}^p_q(w,w)}.
\end{equation}
Moreover, if we assume that 
$w\in A_\infty\cap \mathcal{B}_{p,q}$, 
then for $f\in L^0(\R^n)$, we have that
\begin{equation}\label{160121-9}
\left\|M^{\sharp,d}_{\lambda}f\right\|_{\mathcal{M}^p_q(dx,w)}
\lesssim_{\lambda,p,q,w}
\|f\|_{\mathcal{M}^p_q(dx,w)}.
\end{equation}
\end{proposition}

\begin{proof}
First, we observe that 
for any $\eta>0$,
\begin{equation}\label{151201-1}
M^{\sharp,d}_{\lambda}f(x)
\lesssim_{\lambda,\eta}
M^{(\eta)}f(x).
\end{equation}
To see this, let $Q\in\mathcal{Q}$ and 
$m_f(Q)$ be any median.
We see that by (\ref{151130-1}),
\begin{eqnarray*}
\omega_{\lambda}(f;Q)
&\leq&
\left((f-m_f(Q))\chi_Q\right)^*(\lambda|Q|)\\
&\leq&
(f\cdot\chi_Q)^*\left(\frac{\lambda}{2}|Q|\right)
+
(m_f(Q)\chi_Q)^*\left(\frac{\lambda}{2}|Q|\right)
\leq
2(f\cdot\chi_Q)^*\left(\frac{\lambda}{2}|Q|\right).
\end{eqnarray*}
Moreover, we notice that 
\begin{eqnarray*}
(f\cdot\chi_Q)^*\left(\frac{\lambda}{2}|Q|\right)
\leq
\left(
\frac{2}{\lambda|Q|}\int_0^{\frac{\lambda}{2}|Q|}
\left(f\cdot\chi_Q\right)^*(t)^\eta dt
\right)^\frac{1}{\eta}
\leq
\left(
\frac{2}{\lambda|Q|}\int_Q
|f(x)|^\eta dx
\right)^\frac{1}{\eta}.
\end{eqnarray*}
In summary, we conclude that 
\begin{eqnarray*}
M^{\sharp,d}_{\lambda}f(x)
&\leq&
2\sup_{Q_0\in\mathcal{Q}}\sup_{Q\in\mathcal{D}(Q_0)}
\left(f\cdot\chi_Q\right)^*\left(\frac{\lambda}{2}|Q|\right)\chi_Q(x)\\
&\leq&
2^{1+\frac{1}{\eta}}\lambda^{-\frac{1}{\eta}}
\sup_{Q_0\in\mathcal{Q}}\sup_{Q\in\mathcal{D}(Q_0)}
\left(\frac{1}{|Q|}\int_Q|f(x)|^\eta dx\right)^\frac{1}{\eta}
\chi_Q(x),
\end{eqnarray*}
which implies (\ref{151201-1}).
Next, we choose 
$u\in(1,\infty)$ 
such that 
$w\in A_u$ and 
$\eta:=\frac{q}{u}$.
In view of (\ref{151201-1}) and 
Theorem \ref{th-151202-2}, it follows that
\begin{equation*}
\left\|M^{\sharp,d}_{\lambda}f\right\|_{\mathcal{M}^p_q(w,w)}
\lesssim_{\lambda,\eta}
\left\|M^{(\eta)}f\right\|_{\mathcal{M}^p_q(w,w)}
\lesssim_{p,q,w}
\|f\|_{\mathcal{M}^p_q(w,w)},
\end{equation*}
which proves (\ref{160121-8}).

Meanwhile, if we observe that 
$w\in\mathcal{B}_{p,q}$ is equivalent to 
$w\in\mathcal{B}_{\frac{p}{\eta},\frac{q}{\eta}}$; see \cite{NS15-2}, 
then 
by (\ref{151201-1}) and Theorem $\ref{th-151202-2}$, 
it follows that 
\begin{equation*}
\left\|M^{\sharp,d}_{\lambda}f\right\|_{\mathcal{M}^p_q(dx,w)}
\lesssim_{\lambda,\eta}
\left\|M^{(\eta)}f\right\|_{\mathcal{M}^p_q(dx,w)}
\lesssim_{p,q,w}
\|f\|_{\mathcal{M}^p_q(dx,w)},
\end{equation*}
which proves (\ref{160121-9}).
\end{proof}

By combining 
(\ref{160121-10}),
(\ref{160121-11}),
Corollary \ref{cr-151201-1} and Proposition \ref{pr-151201-2}, 
we obtain Theorem \ref{th-151129-3}.
Next, let us prove Theorem \ref{th-151129-4}.
We will employ a similar method to the proof of Theorem \ref{th-151129-2}.
\begin{proof}[Proof of Theorem \ref{th-151129-4}]
First, let us show that
\begin{equation}\label{151201-6}
\|f\|_{\mathcal{M}^p_s(w,w)}
\lesssim_{p,s,w}
\|M^{\sharp,d}_{\lambda_w}f\|_{\mathcal{M}^p_s(w,w)}
\end{equation}
by assuming $\lim_{l\to\infty}m_f(2^lQ_0)\to0$ for all $Q_0\in\mathcal{Q}$.
To this end, we fix any $Q_0\in\mathcal{Q}$ and
calculate that 
\begin{eqnarray*}
w(Q_0)^{\frac{1}{p}-\frac{1}{s}}
\left(\int_{Q_0}|f(x)|^sdw(x)\right)^\frac{1}{s}
&\lesssim_s&
w(Q_0)^{\frac{1}{p}-\frac{1}{s}}
\left(\int_{Q_0}|f(x)-m_f(2^lQ_0)|^sdw(x)\right)^\frac{1}{s}\\
&&+
w(Q_0)^\frac{1}{p}|m_f(2^lQ_0)|.
\end{eqnarray*} 
By the assumption, it follows that 
\begin{equation}\label{151202-2}
w(Q_0)^{\frac{1}{p}-\frac{1}{s}}
\left(\int_{Q_0}|f(x)|^sdw(x)\right)^\frac{1}{s}
\lesssim_s
\limsup_{l\to\infty}
w(Q_0)^{\frac{1}{p}-\frac{1}{s}}
\left(\int_{Q_0}|f(x)-m_f(2^lQ_0)|^sdw(x)\right)^\frac{1}{s},
\end{equation}
and hence, we focus on the quantity: 
$\left(\int_{Q_0}|f(x)-m_f(2^lQ_0)|^sdw(x)\right)^\frac{1}{s}$.
We employ Proposition \ref{pr-151201-1} again to 
decompose $f-m_f(2^lQ_0)$ and obtain that
\begin{equation}\label{151202-3}
\left(\int_{Q_0}|f(x)-m_f(2^lQ_0)|^sdw(x)\right)^\frac{1}{s}
\lesssim_s
\left(\int_{Q_0}M^{\sharp,d}_{\lambda_w}f(x)^sdw(x)\right)^\frac{1}{s}
+{\rm I},
\end{equation}
where we defined 
\begin{equation*}
{\rm I}
:=
\left[\int_{Q_0}
\left(
\sum_{k\in\N_0}\sum_{j\in J_k}\omega_{\lambda_w}(f;Q^k_j)\chi_{Q^k_j}(x)
\right)^sdw(x)\right]^\frac{1}{s}.
\end{equation*}
Here, we remark that 
the family $\{Q^k_j\}_{k\in\N_0,j\in J_k}\subset \mathcal{D}(2^lQ_0)$ 
is a $w$-sparse family generated by $2^lQ_0$. 
To evaluate ${\rm I}$, 
by considering the suitable dyadic setting,
we have only to calculate the following two terms intrinsically:
\begin{eqnarray*}
{\rm I}_a
&:=&
\left[\int_{Q_0}
\left(
\sum_{k\in\N_0}
\sum_{\substack{j\in J_k:\\ Q^k_j\subset Q_0}}
\omega_{\lambda_w}(f;Q^k_j)\chi_{Q^k_j}(x)
\right)^sdw(x)\right]^\frac{1}{s},\\
{\rm I}_b
&:=&
\left[\int_{Q_0}
\left(
\sum_{k\in\N_0}
\sum_{\substack{j\in J_k:\\ Q^k_j\supsetneq Q_0}}
\omega_{\lambda_w}(f;Q^k_j)\chi_{Q^k_j}(x)
\right)^sdw(x)\right]^\frac{1}{s}.
\end{eqnarray*}
For ${\rm I}_a$, we employ the same method used in 
the proof of Theorem \ref{th-151129-2}, 
namely, by observing 
$\chi_{Q^k_j}\leq C_wM_w[\chi_{Q^k_j\cap \Omega_{k+1}^c}]$
and 
using the boundedness of $M_w$, it follows that 
\begin{equation*}
{\rm I}_a
\lesssim_{s,w}
\left[\int_{Q_0}
\left(
\sum_{k\in\N_0}
\sum_{\substack{j\in J_k:\\ Q^k_j\subset Q_0}}
\omega_{\lambda_w}(f;Q^k_j)\chi_{Q^k_j\cap \Omega_{k+1}^c}(x)
\right)^sdw(x)\right]^\frac{1}{s}.
\end{equation*}
If we notice that 
the summation is taken over the cubes contained in $Q_0$, 
then the disjointness of 
$\{Q^k_j\cap \Omega_{k+1}^c\}_{k\in\N_0,j\in J_k}$ 
yields that 
\begin{equation}\label{151202-4}
{\rm I}_a
\lesssim_{s,w}
\left(\int_{Q_0}M^{\sharp,d}_{\lambda_w}f(x)^sdw(x)\right)^\frac{1}{s}.
\end{equation}
Meanwhile, for ${\rm I}_b$, 
by recalling that $Q^k_j\subset 2^lQ_0$ 
and the dyadic property, 
we can  rewrite the summation of ${\rm I}_b$ as follows:
\begin{equation*}
{\rm I}_b
\leq
\left[\int_{Q_0}
\left(
\sum_{m=1}^{l}
\omega_{\lambda_w}(f;Q_0^{(m)})\chi_{Q_0^{(m)}}(x)
\right)^sdw(x)\right]^\frac{1}{s}.
\end{equation*} 
Here, $Q_0^{(m)}$ denotes the dyadic $m$-th ancestor of $Q_0$, 
that is, 
$Q_0^{(m)}$ is a unique dyadic cube with respect to 
$2^lQ_0$ whose side length is 
$2^m\ell(Q_0)$ and containing $Q_0$. 
With the relation: 
$Q_0\subset Q_0^{(m)}$ in mind, 
we see that 
\begin{eqnarray*}
{\rm I}_b
&\leq&
w(Q_0)^\frac{1}{s}\sum_{m=1}^l\omega_{\lambda_w}(f;Q_0^{(m)})\\
%detail%%%%%%%%%%%%%%%%%%%%%%%%%%%%%%%%%%%%%%%%%%
%&=&
%w(Q_0)^\frac{1}{s}\sum_{m=1}^l\omega_{\lambda_w}(f;Q_0^{(m)})
%\frac{w(Q_0^{(m)})}{w(Q_0^{(m)})}\\
%%%%%%%%%%%%%%%%%%%%%%%%%%%%%%%%%%%%%%%%%%%%%%
&=&
w(Q_0)^\frac{1}{s}\sum_{m=1}^l
\frac{1}{w(Q_0^{(m)})^\frac{1}{s}}
\left(
\int_{Q_0^{(m)}}\omega_{\lambda_w}(f;Q_0^{(m)})^sdw(x)
\right)^\frac{1}{s}\\
&\leq&
w(Q_0)^\frac{1}{s}\sum_{m=1}^l
\frac{1}{w(Q_0^{(m)})^\frac{1}{p}}
\left\|M^{\sharp,d}_{\lambda_w}f\right\|_{\mathcal{M}^p_s(w,w)}.
\end{eqnarray*} 
Now, let us show that 
\begin{equation}\label{151202-1}
\sum_{m=1}^l\frac{1}{w(Q_0^{(m)})^\frac{1}{p}}
\lesssim_{p,w}
\frac{1}{w(Q_0)^\frac{1}{p}}.
\end{equation}
Since $w\in A_\infty$, 
there exist $L_w\in \N$ and 
$\alpha_w>1$ such that 
\begin{equation}\label{160117-4}
w(Q_0^{(i+L_w)})\geq
\alpha_ww(Q_0^{(i)}),
\quad
(i\in\N_0, Q_0\in\mathcal{D}).
\end{equation} 
%For the detail of this property, see Appendix.
In particular, for $i,j\in\N$, it holds that
\begin{equation*}
w(Q_0^{(i+jL_w)})
\geq
\alpha_w
w(Q_0^{(i+(j-1)L_w)})
\geq
\cdots
\geq
\alpha_w^jw(Q_0^{(i)}).
\end{equation*}
Hence, the left-hand side of (\ref{151202-1}) can be controlled as follows:
\begin{eqnarray*}
\sum_{m=1}^l\frac{1}{w(Q_0^{(m)})^\frac{1}{p}}
\leq
\sum_{i=1}^{L_w}
\sum_{j=0}^{[l/L_w]+1}
\frac{1}{w(Q_0^{(i+jL_w)})^\frac{1}{p}}
\leq
\sum_{i=1}^{L_w}
\frac{1}{w(Q_0^{(i)})^\frac{1}{p}}
\sum_{j=0}^{[l/L_w]+1}
\frac{1}{\alpha_w^{\frac{j}{p}}}.
\end{eqnarray*}
If we recall that 
$w(Q_0)\leq w(Q_0^{(i)})$ for any $i\in\N$ and
that $\alpha_w>1$, 
then we see that 
\begin{equation*}
\sum_{m=1}^l\frac{1}{w(Q_0^{(m)})^\frac{1}{p}}
\lesssim_{p,w}
\frac{L_w}{w(Q_0)^\frac{1}{p}},
\end{equation*}
which implies (\ref{151202-1}), since 
the constant $L_w$ and $\alpha_w$ depend only on $w$ (and dimension $n$).
Remark that the implicit constant in (\ref{151202-1}) 
does not depend on $l\in\N$ as well.
Therefore, we obtain the estimate of ${\rm I}_b$:
\begin{equation}\label{151202-5}
{\rm I}_b
\lesssim_{p,w}
w(Q_0)^{-\frac{1}{p}+\frac{1}{s}}
\left\|M^{\sharp,d}_{\lambda_w}f\right\|_{\mathcal{M}^p_s(w,w)}.
\end{equation}
In total, 
it follows from (\ref{151202-3}), (\ref{151202-4}) and (\ref{151202-5})
that 
\begin{eqnarray*}
\left(\int_{Q_0}|f(x)-m_f(2^lQ_0)|^sdw(x)\right)^\frac{1}{s}
&\lesssim_{p,s,w}&
\left(\int_{Q_0}M^{\sharp,d}_{\lambda_w}f(x)^sdw(x)\right)^\frac{1}{s}\\
&&+
w(Q_0)^{-\frac{1}{p}+\frac{1}{s}}
\left\|M^{\sharp,d}_{\lambda_w}f\right\|_{\mathcal{M}^p_s(w,w)}
\end{eqnarray*}
and hence, in view of (\ref{151202-2}), we conclude (\ref{151201-6}).

Next, assuming the weighted integral condition (\ref{weight-integral})
for $p,q$ and $w$, we prove 
\begin{equation}\label{151202-6}
\|f\|_{\mathcal{M}^p_q(dx,w)}
\lesssim_{p,q,w}
\left\|M^{\sharp,d}_{\lambda_w}f\right\|_{\mathcal{M}^p_q(dx,w)}.
\end{equation}
Another inequality:
\begin{equation*}
\sup_{Q\in\mathcal{Q}}
\Phi_{p,q,w}(Q)
\left(\frac{1}{w(Q)}\int_Q|f(x)|^sdw(x)\right)^\frac{1}{s}
\leq
\|f\|_{\mathcal{M}^p_q(dx,w)}
\end{equation*}
follows from H\"{o}lder's inequality directly.
We need only mimic the above proof.
Fix any $Q_0\in\mathcal{D}$ and calculate that 
\begin{equation*}
|Q_0|^{\frac{1}{p}-\frac{1}{q}}
\left(\int_{Q_0}|f(x)|^qdw(x)\right)^\frac{1}{q}
\lesssim_q
\limsup_{l\to\infty}
|Q_0|^{\frac{1}{p}-\frac{1}{q}}
\left(\int_{Q_0}|f(x)-m_f(2^lQ_0)|^qdw(x)\right)^\frac{1}{q}.
\end{equation*}
Again, Proposition \ref{pr-151201-1}, 
Lerner's decomposition formula
reduces 
the matters to show the inequality:
$
{\rm II}_a+{\rm II}_b
\lesssim_{p,q,w}
\|M^{\sharp,d}_{\lambda_w}f\|_{\mathcal{M}^p_q(dx,w)},
$
where,
\begin{eqnarray*}
{\rm II}_a
&:=&
|Q_0|^{\frac{1}{p}-\frac{1}{q}}
\left[\int_{Q_0}
\left(
\sum_{k\in\N_0}
\sum_{\substack{j\in J_k:\\ Q^k_j\subset Q_0}}
\omega_{\lambda_w}(f;Q^k_j)\chi_{Q^k_j}(x)
\right)^sdw(x)\right]^\frac{1}{s},\\
{\rm II}_b
&:=&
|Q_0|^{\frac{1}{p}-\frac{1}{q}}
\left[\int_{Q_0}
\left(
\sum_{k\in\N_0}
\sum_{\substack{j\in J_k:\\ Q^k_j\supsetneq Q_0}}
\omega_{\lambda_w}(f;Q^k_j)\chi_{Q^k_j}(x)
\right)^sdw(x)\right]^\frac{1}{s},
\end{eqnarray*}
and $\{Q^k_j\}_{k\in\N_0,j\in J_k}\subset\mathcal{D}(2^lQ_0)$ 
is a $w$-sparse family generated by 
$2^lQ_0$.
For ${\rm II}_a$, 
we may apply the same argument used to 
evaluate ${\rm I}_a$; see (\ref{151202-4})
to obtain that 
\begin{equation*}
{\rm II}_a
\lesssim_{q,w}
|Q_0|^{\frac{1}{p}-\frac{1}{q}}
\left(\int_{Q_0}M^{\sharp,d}_{\lambda_w}f(x)^qdw(x)
\right)^\frac{1}{q}
\leq
\left\|M^{\sharp,d}_{\lambda_w}f\right\|_{\mathcal{M}^p_q(dx,w)}.
\end{equation*}
For ${\rm II}_b$, we employ the weighted integral condition 
(\ref{weight-integral})
for $p,q$ and $w$ instead of the $A_\infty$ condition.
As before, using the dyadic $m$-ancestor of $Q_0$, 
we see that 
\begin{eqnarray*}
{\rm II}_b
&\leq&
|Q_0|^{\frac{1}{p}-\frac{1}{q}}
\left(
\int_{Q_0}
\left(
\sum_{m=1}^l\omega_{\lambda_w}(f;Q_0^{(m)})\chi_{Q_0^{(m)}}(x)
\right)^qdw(x)
\right)^\frac{1}{q}\\
&=&
\Phi_{p,q,w}(Q_0)
\sum_{m=1}^l\omega_{\lambda_w}(f;Q_0^{(m)})\\
%detai%%%%%%%%%%%%%%%%%%%%%%%%%%%%%%%%%%%%%%
%&=&
%\Phi_{p,q,w}(Q_0)
%\sum_{m=1}^l
%\frac{1}{w(Q_0^{(m)})^\frac{1}{q}}
%\left(
%\int_{Q_0^{(m)}}\omega_{\lambda_w}(f;Q_0^{(m)})^qdw(x)
%\right)^\frac{1}{q}\\
%%%%%%%%%%%%%%%%%%%%%%%%%%%%%%%%%%%%%%%%%%
&\leq&
\Phi_{p,q,w}(Q_0)
\sum_{m=1}^l
\frac{1}{w(Q_0^{(m)})^\frac{1}{q}}
\left(
\int_{Q_0^{(m)}}M^{\sharp,d}_{\lambda_w}f(x)^qdw(x)
\right)^\frac{1}{q}.
\end{eqnarray*}
From the definition of the norm of $\mathcal{M}^p_q(dx,w)$, 
and the weighted integral condition
(\ref{weight-integral})
it follows that
\begin{eqnarray*}
{\rm II}_b
\leq
\Phi_{p,q,w}(Q_0)
\sum_{m=1}^l\frac{1}{\Phi_{p,q,w}(Q_0^{(m)})}
\left\|M^{\sharp,d}_{\lambda}f\right\|_{\mathcal{M}^p_q(dx,w)}
\lesssim_{p,q,w}
\left\|M^{\sharp,d}_{\lambda}f\right\|_{\mathcal{M}^p_q(dx,w)}.
\end{eqnarray*}
\end{proof}

By invoking 
Lemma \ref{lm-151201-3} in Chapter 2 and (\ref{151201-4}), 
we obtain the direct analogy from 
original type of the sharp maximal inequality :
\begin{corollary}\label{cr-151204-1}
Let $0<q\leq p<\infty$ and $w\in A_\infty$.
\begin{enumerate}
\item
(Komori-Shirai type)
For any $f\in L^0(\R^n)$ satisfying 
$Mf\in\mathcal{M}^{p_0}_{q_0}(w_0,w_0)$ 
for some $0<q_0\leq p_0<\infty$ and $w_0\in A_\infty$, 
we have that
\begin{equation*}
\|f\|_{\mathcal{M}^p_q(w,w)}
\sim_{p,q,w}
\left\|M^{\sharp,d}_{\lambda_w}f\right\|_{\mathcal{M}^p_q(w,w)}
\lesssim_{\lambda_w}
\left\|f^\sharp\right\|_{\mathcal{M}^p_q(w,w)}.
\end{equation*}
\item
(Samko type)
For any $f\in L^0(\R^n)$ satisfying 
$Mf\in \mathcal{M}^{p_0}_{q_0}(dx,w_0)$ for some 
$0<q_0\leq p_0<\infty$ and $w_0$ satisfying 
the weighted integral condition (\ref{weight-integral})
for $p_0,q_0$ and $w_0$, we have that
\begin{equation*}
\|f\|_{\mathcal{M}^p_q(dx,w)}
\sim_{p,q,w}
\left\|M^{\sharp,d}_{\lambda_w}f\right\|_{\mathcal{M}^p_q(dx,w)}
\lesssim_{\lambda_w}
\left\|f^\sharp\right\|_{\mathcal{M}^p_q(dx,w)}.
\end{equation*}
%Assume the weighted integral condition for $p,q$ and $w$.
%Then the inequality 
%\begin{equation*}
%\|f\|_{\mathcal{M}^p_q(dx,w)}
%\sim_{p,q,w}
%\|M^{\sharp,d}_{\lambda_w}f\|_{\mathcal{M}^p_q(dx,w)}
%\lesssim_{\lambda_w}
%\|M^{\sharp}_{(1)}f\|_{\mathcal{M}^p_q(dx,w)}
%\end{equation*}
%holds for any $f\in L^0(\R^n)$ 
%satisfying $Mf\in\mathcal{M}^{p_0}_{q_0}(dx,w_0)$ 
%for some $0<q_0\leq p_0<\infty$ and $w_0\in A_\infty$
%with the weak weighted integral condition for $p_0,q_0$ and $w_0$: (\ref{151203-2}).
\end{enumerate}
\end{corollary}

Let us compare the above results 
with some recent researches.
Sawano and Tanaka in \cite{SaTa07-2} proved the following refinement of
sharp maximal inequality:
\begin{theorem}{\rm \cite[Theorem 1.3]{SaTa07-2}}\label{th-151201-4}
Let $1<q\leq p<\infty$. Then for any $f\in L^0(\R^n)$, 
we have that
\begin{equation*}
\|Mf\|_{\mathcal{M}^p_q(dx,dx)}
\lesssim_{p,q}
\left\|f^\sharp\right\|_{\mathcal{M}^p_q(dx,dx)}
+
\|f\|_{\mathcal{M}^p_1(dx,dx)}.
\end{equation*}
\end{theorem}
In view of (\ref{151201-4}) and Theorem \ref{th-151202-2}  
Theorem \ref{th-151129-3} improves and generalizes Theorem \ref{th-151201-4}.
Meanwhile, there exists a weighted result for the sharp maximal inequality 
obtained by Komori-Furuya in \cite{Komori15}:
\begin{theorem}{\rm \cite[Theorem 12]{Komori15}}\label{th-151201-5}
Let $1<q_0\leq q\leq p<\infty$ and $w\in A_\infty$.
Then for $f\in L^0(\R^n)$ satisfying $Mf\in \mathcal{M}^p_{q_0}(w,w)$, 
it holds that 
\begin{equation*}
\|f\|_{\mathcal{M}^p_q(w,w)}
\lesssim_{p,q,w}
\left\|f^\sharp\right\|_{\mathcal{M}^p_q(w,w)}.
\end{equation*}
\end{theorem}
We can see easily that 
Corollaries \ref{cr-151201-2} and \ref{cr-151204-1} 
improve and generalize Theorem \ref{th-151201-5}.

\section{Proof of Theorems \ref{th-160626-1} and \ref{th-160113-1} including another application}
\label{s5}

%%%%%%%%%%%%%%%%%%%%%%%%%%%%%%%%%%%%%%%%%%%%%%%%
%the definition of the singular integral operator
We will show Theorem \ref{th-160626-1} with more general settings. 
Recall the definition of the singular integral operators.  
A singular integral operator $T$ is a bounded linear operator on the unweighted $L^2$ space, 
for which there exists a kernel $K$ on $\R^n\times\R^n$ satisfying the following conditions:
\begin{enumerate}
\item
(Size condition)
There exists $C>0$ such that $|K(x,y)|\leq \frac{C}{|x-y|^n}$ for any $x\neq y$.
\item
(H\"{o}rmander condition)
For some $\theta\in(0,1]$,
\begin{equation*}
|K(x+h,y)-K(x,y)|+|K(x,y+h)-K(x,y)|\leq \frac{C|h|^\theta}{|x-y|^{n+\theta}},
\quad \mbox{if}\  |x-y|>2|h|.
\end{equation*}
\item
If $f\in L^\infty_{\rm c}$, the set of all compactly supported
$L^\infty$-functions, then
\begin{equation*}
Tf(x)=\int_{\mathbb{R}^n}K(x,y)f(y)dy
\quad
(x\notin \mbox{supp}(f)).
\end{equation*}
\end{enumerate}
As Coifman and Fefferman showed 
in \cite[Theorems I and III]{CF74}, 
$T$ is well defined 
as a bounded operator on $L^q(w)$, when $w\in A_q$ with $1<q<\infty$.
%%%%%%%%%%%%%%%%%%%%%%%%%%%%%%%%%%%%%%%%%%%%%%

We need to recall some examples 
of the singular integral operators.
\begin{definition}[\cite{BGST10}]
A singular integral operator $T$ is called 
a genuine singular integral operator 
if there exist constants $C,\theta>0$ and $R\in O(n)$ such that
for any $x,y\in\R^n$ satisfying 
$x-y\in V_{\theta,R}:=R\{u=(u',u_n)\in\R^n:|u'|<\theta|u_n|\}$, 
the integral kernel $K$ satisfies
\begin{equation*}
K(x,y)\geq
\frac{C}{|x-y|^n},
\end{equation*}
where $O(n)$ denotes the set of all orthogonal matrices in $\R^n$.
\end{definition}
A typical example of the genuine singular integral operators is 
the Riesz transform. 
Our aim in this section is to show the following theorem:
\begin{theorem}\label{th-160104-1}
Let $1<q\leq p<\infty$ and $w$ be a weight.
\begin{enumerate}
\item
Assume that 
$M$ is bounded on $\mathcal{M}^p_q(dx,w)$ 
and the weighted integral condition 
(\ref{weight-integral})
for $p,q$ and $w$.
Then 
we can extend any 
singular integral operator 
$T_0$ initially defined on 
$L^\infty_{\rm c}$
to a bounded linear operator $T$ on $\mathcal{M}^p_q(dx,w)$:
\begin{equation*}
\|Tf\|_{\mathcal{M}^p_q(dx,w)}
\leq
C\|f\|_{\mathcal{M}^p_q(dx,w)}
\quad
(f\in\mathcal{M}^p_q(dx,w)).
\end{equation*}
\item
As a converse assertion, we have the following:
\begin{enumerate}
\item
If the Riesz transform 
$R_i$ is bounded on 
$\mathcal{M}^p_q(dx,w)$ 
for some $i=1,\ldots,n$,
then 
$w$ satisfies (\ref{151221-1}).
\item
Assume that $w\in\mathcal{B}_{p,q}$ satisfies the doubling condition.
If there exists a genuine singular integral operator $T$ which is 
bounded on 
$\mathcal{M}^p_q(dx,w)$, 
then the weighted integral condition (\ref{weight-integral}) 
for $p,q$ and $w$ holds.
\end{enumerate}
\end{enumerate} 
\end{theorem}

Once we prove Theorem \ref{th-160104-1}, 
then we obtain Theorem \ref{th-160626-1} as a corollary. 
To see this, we admit Theorem \ref{th-160104-1} for a moment. 
If the Riesz transform $R_i$ is bounded on $\mathcal{M}^p_q(dx,w)$, 
then 
\eqref{151221-1} holds by the assertion (2)-(a) in Theorem \ref{th-160104-1}. 
Then by Lemma \ref{lm-160113-1} or Remark \ref{rm-160628-1}, 
we see that $w\in \mathcal{B}_{p,q}\cap A_{q+1}$, particularly, 
$w$ is doubling. 
Here, we recall that the Riesz transform $R_i$ 
is an example of the genuin singular integral operators. 
Hence, we may apply the assertion (2)-(b) in Theorem \eqref{th-160104-1} to 
obtain the weighted integral condition \eqref{weight-integral} for $p,q$ and $w$. 
This means the assertion \eqref{item-160626-2} in Theorem \ref{th-160626-1}. 
The assertion \eqref{item-160626-1} in Theorem \ref{th-160626-1} is clear from  
the assertion (1) in Theorem \ref{th-160104-1}. 

Since the proof of the assertion (1) of Theorem \ref{th-160104-1} 
is easier than the one of Theorem \ref{th-160113-1}, 
we prove only the assertion (2) of Theorem \ref{th-160104-1} here.
\begin{proof}[Proof of Theorem \ref{th-160104-1}-$(2)$-$(a)$]
Assume that 
$R_i$ is bounded on $\mathcal{M}^p_q(dx,w)$ for some $i=1,\ldots,n$.
To show (\ref{151221-1}), we fix any $Q\in\mathcal{Q}$ 
and take any $f\geq0$ such that $f\cdot\chi_Q\in\mathcal{M}^p_q(dx,w)$.
We focus on the cube $\widetilde{Q}:=Q+2\ell(Q)e_i$. 
Here, $e_i$ denotes the $i$-th elementary vector.
If we notice that for any $x\in\widetilde{Q}$, 
it holds that 
$x_i-y_i\geq\ell(Q)$ and that 
$|x-y|\leq C\ell(Q)$, 
then
we have the pointwise estimate
\begin{equation*}
\left|R_i[f\cdot\chi_Q](x)\right|
=
\int_Qf(y)\frac{x_i-y_i}{|x-y|^{n+1}}dy
\geq
C\frac{1}{|Q|}\int_Qf(y)dy.
\end{equation*}
Hence
it follows from  the boundedness of $R_i$ that 
\begin{equation*}
\frac{1}{|Q|}\int_Qf(y)dy
%\leq
%C\inf_{x\in\widetilde{Q}}|R_i[f\cdot\chi_Q](x)|
\leq
\frac{C}{\|\chi_{\widetilde{Q}}\|_{\mathcal{M}^p_q(dx,w)}}
\left\|R_i[f\cdot\chi_Q]\chi_{\widetilde{Q}}\right\|_{\mathcal{M}^p_q(dx,w)}
\leq
C
\frac{\|f\cdot\chi_Q\|_{\mathcal{M}^p_q(dx,w)}}
{\|\chi_{\widetilde{Q}}\|_{\mathcal{M}^p_q(dx,w)}}.
\end{equation*}
Namely, we have that 
for any $f\geq0$ such that 
$f\cdot\chi_Q\in\mathcal{M}^p_q(dx,w)$, 
\begin{equation}\label{160114-7}
\frac{1}{|Q|}\int_Qf(y)dy
\cdot
\|\chi_{\widetilde{Q}}\|_{\mathcal{M}^p_q(dx,w)}
\leq
C
\|f\cdot\chi_Q\|_{\mathcal{M}^p_q(dx,w)}.
\end{equation}
Now, we take any $g\geq0$ such that 
$\|g\cdot\chi_{\widetilde{Q}}\|_{\mathcal{M}^p_q(dx,w)}\leq1$ 
and 
put $f=\left|R_i[g\cdot\chi_{\widetilde{Q}}]\right|$ as in (\ref{160114-7})
to obtain that
\begin{equation}\label{160114-8}
\frac{1}{|Q|}\int_Q\left|R_i[g\cdot\chi_{\widetilde{Q}}](x)\right|dx
\cdot
\|\chi_{\widetilde{Q}}\|_{\mathcal{M}^p_q(dx,w)}
\leq
C
\|g\cdot\chi_{\widetilde{Q}}\|_{\mathcal{M}^p_q(dx,w)}
\leq
C.
\end{equation}
Here, we again used the boundedness of $R_i$.
Meanwhile, 
if we go through a similar argument as before, 
we notice
\begin{equation}\label{160114-9}
\frac{1}{|\widetilde{Q}|}\int_{\widetilde{Q}}g(y)dy
\leq
C\inf_{x\in Q}\left|R_i[g\cdot\chi_{\widetilde{Q}}](x)\right|
\leq
C\frac{1}{|Q|}\int_Q\left|R_i[g\cdot\chi_{\widetilde{Q}}](x)\right|dx.
\end{equation}
By combining (\ref{160114-8}) and (\ref{160114-9}),  
it follows that for any $g\geq0$ such that 
$\|g\cdot\chi_{\widetilde{Q}}\|_{\mathcal{M}^p_q(dx,w)}\leq1$, 
\begin{equation}\label{160114-10}
\frac{1}{|\widetilde{Q}|}\int_{\widetilde{Q}}g(y)dy
\cdot
\|\chi_{\widetilde{Q}}\|_{\mathcal{M}^p_q(dx,w)}
\leq
C.
\end{equation}
Moreover using (\ref{dual2}), 
we see that
\begin{eqnarray*}
\frac{1}{|\widetilde{Q}|}\|\chi_{\widetilde{Q}}\|_{\mathcal{M}^p_q(dx,w)}
\left\|w^{-\frac{1}{q}}\chi_{\widetilde{Q}}\right\|_{H^{q',n(1-q/p)}}
=
\frac{
\|\chi_{\widetilde{Q}}\|_{\mathcal{M}^p_q(dx,w)}
}{|\widetilde{Q}|}
\sup_{\substack {g\geq0:\\ \|g\chi_{\widetilde{Q}}\|_{\mathcal{M}^p_q(dx,w)}\leq1}}
\int_{\widetilde{Q}}g(y)dy
\leq
C,
\end{eqnarray*}
which implies (\ref{151221-1}).
\end{proof}

Next, we prove the assertion (2)-(b) in Theorem \ref{th-160104-1}. 
That is, we show 
the weighted integral condition
(\ref{weight-integral}) for $p,q$ and $w$
by assuming that $w\in \mathcal{B}_{p,q}$ satisfies the doubling condition 
and that
a certain genuine singular integral operator $T$ is bounded
on $\mathcal{M}^p_q(dx,w)$.
To this end, we recall the equivalent condition 
of the weighted integral condition (\ref{weight-integral}).

\begin{lemma}[\cite{NS15}]\label{lm-151222-1}
Let $0<q\leq p<\infty$ and $w\in \mathcal{B}_{p,q}$ satisfies the doubling condition.
Then the following are equivalent:
\begin{enumerate}
\item
The weighted integral condition for $p,q$ and $w$ holds.
\item
There exists $c>1$ such that 
\begin{equation*}
2\Phi_{p,q,w}(Q)
\leq
\Phi_{p,q,w}(cQ)
\end{equation*}
for any $Q\in\mathcal{Q}$.
\end{enumerate} 
\end{lemma}
The idea of the following proof goes back to \cite[Theorem 6.9]{HNS15}.
\begin{proof}[Proof of Theorem \ref{th-160104-1}-$(2)$-$(b)$]
By assuming that 
the weighted integral condition (\ref{weight-integral})
for $p,q$ and $w$ 
fails to hold, 
let us obtain the contradiction.
With Lemma \ref{lm-151222-1} in mind, 
we assume that 
for any $m\in\N$, there exists $Q_m\in\mathcal{Q}$ such 
that $2\Phi_{p,q,w}(Q_m)>\Phi_{p,q,w}(2^mQ_m)$.
We take a genuine singular integral operator $T$ 
and denote the cone $V_{\theta,R}$ by $V$.
Set 
\begin{equation*}
f_m(y)
:=
\chi_{V-c(Q_m)}(-y)\chi_{2^{m-1}Q_m\setminus 2^6Q_m}(y),
\end{equation*}
for $m\geq7$.
Since $V$ is a cone, 
we may take another cube $R_m$ so that 
$
R_m\subset 
\left(V+c(Q_m)\right)\cap \left(2^4Q_m\setminus Q_m\right)$, 
$Q_m\subset C_V R_m$ for some $C_V>0$
and that 
$|R_m|\sim_V|Q_m|$.
Note that for any $x\in R_m$ 
and any $y\in2^{m-1}Q_m\setminus2^6Q_m$ such that
$-y\in V-c(Q_m)$, we have 
$x-y\in V$ and 
$\ell(Q_m)\lesssim 
|x-y|
\lesssim
2^m\ell(Q_m)$.
Hence, it follows that for any 
$x\in R_m$,
\begin{eqnarray*}
Tf_m(x)
&=&
\int_{\R^n}K(x,y)
\chi_{V-c(Q_m)}(-y)\chi_{2^{m-1}Q_m\setminus 2^6Q_m}(y)dy\\
&\sim&
\int_{
\left\{x-y\in V:
\ell(Q_m)<|x-y|\leq2^m\ell(Q_m)
\right\}}
K(x,y)f_m(y)dy\\
&\gtrsim&
\int_{\left\{\ell(Q_m)<|x-y|\leq2^m\ell(Q_m)\right\}}
\frac{dy}{|x-y|^n}\\
&\sim&
m.
\end{eqnarray*}
This implies that 
$
m\cdot\chi_{R_m}(x)
\lesssim
Tf_m(x).
$
Meanwhile, if we recall that 
$Q_m\subset C_V R_m$, 
$|Q_m|\sim|R_m|$ and $w$ is a doubling weight, 
then we see that 
\begin{equation*}
\Phi_{p,q,w}(Q_m)
\lesssim
|R_m|^\frac{1}{p}
\left(\frac{w(C_VR_m)}{|R_m|}\right)^\frac{1}{q}
\lesssim
\Phi_{p,q,w}(R_m)
\leq
\|\chi_{R_m}\|_{\mathcal{M}^p_q(dx,w)}.
\end{equation*}
Hence, the boundedness of $T$ yields 
\begin{equation*}
\Phi_{p,q,w}(Q_m)
\lesssim
\frac{1}{m}\|Tf_m\|_{\mathcal{M}^p_q(dx,w)}
\lesssim
\frac{1}{m}\|\chi_{2^mQ_m}\|_{\mathcal{M}^p_q(dx,w)}.
\end{equation*}
Moreover, by 
$w\in\mathcal{B}_{p,q}$ and 
$\Phi_{p,q,w}(2^mQ_m)<2\Phi_{p,q,w}(Q_m)$, 
we see that
\begin{equation*}
\Phi_{p,q,w}(Q_m)
\lesssim
\frac{1}{m}\Phi_{p,q,w}(2^mQ_m)
\leq
\frac{2}{m}\Phi_{p,q,w}(Q_m),
\end{equation*}
which implies $1\lesssim m^{-1}$ holds for any $m\in\N$.
This is a contradiction.
\end{proof}

Next, 
we turn into the proof of Theorem \ref{th-160113-1}.
Again, we prove Theorem \ref{th-160113-1} with more general setting as follows. 

\begin{theorem}\label{th-160626-2}
Let $1<q\leq p<\infty$, $b\in {\rm BMO}$ and $T$ 
be a singular integral operator.  
Assume that $M$ is bounded on $\mathcal{M}^p_q(dx,w)$ 
and the weighted integral condition \eqref{weight-integral} for $p,q$ and $w$. 
Then we can extend the commutator $[b,T]$ as a bounded linear operator on 
$\mathcal{M}^p_q(dx,w)$: 
\[
\left\|
[b,T]f
\right\|_{\mathcal{M}^p_q(dx,w)}
\leq
C
\|f\|_{\mathcal{M}^p_q(dx,w)}
\quad
(f\in \mathcal{M}^p_q(dx,w)).
\]
\end{theorem}

To show Theorem \ref{th-160626-2}, we first extend the definition 
of the commutator $[b,T]$ as a linear operator 
defined on $\mathcal{M}^p_q(dx,w)$.
We remark that 
once we assume the boundedness of $M$ on $\mathcal{M}^p_q(dx,w)$, 
then 
$w\in A_\infty$ automatically by Lemma \ref{lm-160113-1}.
Let $b\in {\rm BMO}$ and $T$ be a singular integral operator.
Moreover, we take $\varepsilon\in(0,1)$ as in Corollary \ref{cr-160113-4}.
For $f\in\mathcal{M}^p_q(dx,w)$ and $x\in\R^n$, 
we define 
\begin{equation*}
[b,T]f(x)
:=
[b,T]_0\left(f\cdot\chi_{2Q}\right)(x)
+
\int_{\R^n\setminus 2Q}
\left(b(x)-b(y)\right)K(x,y)f(y)dy,
\end{equation*}
where $Q$ is any cube containing the point $x$ 
and 
$[b,T]_0$ denotes the commutator as a bounded 
linear operator on $L^\frac{1}{1-\varepsilon}(\R^n)$.
It is easy to check that 
the definition of $[b,T]f(x)$ does not depend on 
the choice of the cube $Q$ which contains $x$.
In addition, since we know that 
$f\cdot\chi_Q\in L^\frac{1}{1-\varepsilon}(\R^n)$
for any $f\in\mathcal{M}^p_q(dx,w)$ and $Q\in\mathcal{Q}$
by Corollary \ref{cr-160113-4}, the first term: 
$[b,T]_0\left(f\cdot\chi_{2Q}\right)(x)$ is well defined. 
For the second term, 
we have the following lemma:
\begin{lemma}\label{lm-160114-9}
Let $1<q\leq p<\infty$ and assume the boundedness of $M$ on $\mathcal{M}^p_q(dx,w)$ 
and  
the weighted integral condition (\ref{weight-integral})
for $p,q$ and $w$. 
Then for any $x\in\R^n$ and any $Q\in\mathcal{Q}$ containing $x$, 
we have that
\begin{equation*}
\int_{\R^n\setminus 2Q}
\left|\left(b(x)-b(y)\right)K(x,y)f(y)\right|
dy
\leq
C
\frac{\|f\|_{\mathcal{M}^p_q(dx,w)}}{\Phi_{p,q,w}(Q)}
\left(
\|b\|_{\rm BMO}
+
|b(x)-b_Q|
\right),
\end{equation*}
where $b_Q$ denotes the mean value of $b$ 
on $Q$ with respect to the Lebesgue measure, 
that is, 
$b_Q:=\frac{1}{|Q|}\int_Qb(x)dx$.
\end{lemma}

\begin{proof}
It follows from the size condition of the kernel $K$ that
\begin{equation*}
\int_{\R^n\setminus 2Q}
\left|\left(b(x)-b(y)\right)K(x,y)f(y)\right|
dy
\leq
C({\rm I}+{\rm II}),
\end{equation*}
where we denined
\begin{equation*}
{\rm I}
:=
\int_{\R^n\setminus 2Q}\frac{|b(x)-b_Q|}{|x-y|^n}|f(y)|dy,
\quad
{\rm II}
:=
\int_{\R^n\setminus 2Q}\frac{|b(y)-b_Q|}{|x-y|^n}|f(y)|dy.
\end{equation*}
We first evaluate the second term {\rm II}.
Since we assume the boundedness of $M$ on 
$\mathcal{M}^p_q(dx,w)$, 
we can choose $\varepsilon\in(0,1)$ such that 
$M^{\left(\frac{1}{1-\varepsilon}\right)}$ 
is also bounded on $\mathcal{M}^p_q(dx,w)$ 
by Lemma \ref{lm-160113-2}.
Then we have 
\begin{eqnarray*}
{\rm II}
&\leq&
\sum_{l=1}^\infty
\frac{C}{|2^lQ|}
\int_{2^lQ}|b(y)-b_Q||f(y)|dy\\
&\leq&
C
\sum_{l=1}^\infty
\left(
\frac{1}{|2^lQ|}
\int_{2^lQ}|b(y)-b_Q|^{\frac{1}{\varepsilon}}dy
\right)^\varepsilon
\left(\frac{1}{|2^lQ|}\int_{2^lQ}|f(y)|^{\frac{1}{1-\varepsilon}}dy\right)^{1-\varepsilon}\\
&\leq&
C
\sum_{l=1}^\infty
\left(
\frac{1}{|2^lQ|}
\int_{2^lQ}|b(y)-b_Q|^{\frac{1}{\varepsilon}}dy
\right)^\varepsilon
\frac{\|\chi_{2^lQ}\|_{\mathcal{M}^p_q(dx,w)}}
{\Phi_{p,q,w}(2^lQ)}
\inf_{x\in2^lQ}M^{\left(\frac{1}{1-\varepsilon}\right)}f(x)\\
&\leq&
C
\sum_{l=1}^\infty
\left(
\frac{1}{|2^lQ|}
\int_{2^lQ}|b(y)-b_Q|^{\frac{1}{\varepsilon}}dy
\right)^\varepsilon
\frac{\|f\|_{\mathcal{M}^p_q(dx,w)}}
{\Phi_{p,q,w}(2^lQ)}.
\end{eqnarray*}
Now, we focus on the term related to the function $b$. 
It is easy to see that 
$|b_{2^lQ}-b_Q|\leq
Cl\|b\|_{\rm BMO}$. 
Hence, we have that 
\begin{eqnarray*}
\left(
\frac{1}{|2^lQ|}
\int_{2^lQ}|b(y)-b_Q|^{\frac{1}{\varepsilon}}dy
\right)^\varepsilon
&\leq&
\left(
\frac{1}{|2^lQ|}
\int_{2^lQ}|b(y)-b_{2^lQ}|^{\frac{1}{\varepsilon}}dy
\right)^\varepsilon
+
|b_{2^lQ}-b_Q|\\
&\leq&
C(1+l)\|b\|_{\rm BMO}
\end{eqnarray*}
for all $l\in\N$.
It thus follows from Lemma \ref{lm-160115-1}
that 
\begin{equation*}
{\rm II}
\leq
C\|b\|_{\rm BMO}
\|f\|_{\mathcal{M}^p_q(dx,w)}
\sum_{l=1}^\infty
\frac{l}{\Phi_{p,q,w}(2^lQ)}
\leq
C\|b\|_{\rm BMO}
\frac{\|f\|_{\mathcal{M}^p_q(dx,w)}}
{\Phi_{p,q,w}(Q)}.
\end{equation*}

Next, we focus on the first term ${\rm I}$.
We employ 
(\ref{dual}),
(\ref{151221-1}) 
and the weighted integral condition (\ref{weight-integral})
for $p,q$ and $w$
to obtain that
\begin{eqnarray*}
{\rm I}
&\leq&
C
|b(x)-b_Q|
\left(
\sum_{l=1}^\infty
\frac{1}{|2^lQ|}
\int_{2^lQ}|f(y)|dy
\right)\\
&\leq&
C
|b(x)-b_Q|
\left(
\sum_{l=1}^\infty
\frac{1}{|2^lQ|}
\left\|w^{-\frac{1}{q}}\chi_{2^lQ}\right\|_{H^{q',n(1-q/p)}}
\|f\|_{\mathcal{M}^p_q(dx,w)}
\right)\\
&\leq&
C
|b(x)-b_Q|
\cdot
\|f\|_{\mathcal{M}^p_q(dx,w)}
\sum_{l=1}^\infty\frac{1}{\Phi_{p,q,w}(2^lQ)}\\
&\leq&
C
|b(x)-b_Q|
\frac{\|f\|_{\mathcal{M}^p_q(dx,w)}}{\Phi_{p,q,w}(Q)}.
\end{eqnarray*}
Hence, we complete the proof of Lemma \ref{lm-160114-9}.
\end{proof}

Additionally, we need the following simple observation.
\begin{lemma}\label{lm-160115-4}
Let $0<q<\infty$, $w\in A_\infty$ and $b\in {\rm BMO}$.
Then for any cube $Q$, we have that 
\begin{equation*}
\left(
\frac{1}{w(Q)}
\int_Q
|b(x)-b_Q|^qw(x)dx
\right)^\frac{1}{q}
\leq
C
\|b\|_{\rm BMO}.
\end{equation*}
\end{lemma}

For the proof of Lemma \ref{lm-160115-4}, 
see \cite{grafakos-book} for example. 
Using Theorem \ref{th-151129-4}
and 
Lemmas \ref{lm-160114-9} and 
\ref{lm-160115-4}, 
let us complete the proof of Theorem \ref{th-160626-2}.

\begin{proof}[Proof of Theorem \ref{th-160626-2}]
Fix any $f\in\mathcal{M}^p_q(dx,w)$ and
any cube $Q$ and set
\begin{eqnarray*}
{\rm I}
&:=&
|Q|^{\frac{1}{p}-\frac{1}{q}}
\left(\int_Q|[b,T]_0(f\cdot\chi_{2Q})(x)|^q w(x)dx
\right)^\frac{1}{q},\\
{\rm II}
&:=&
|Q|^{\frac{1}{p}-\frac{1}{q}}
\left(\int_Q
\left|
\int_{\R^n\setminus 2Q}(b(x)-b(y))K(x,y)f(y)dy
\right|^q w(x)dx
\right)^\frac{1}{q}.
\end{eqnarray*}
Here, 
$[b,T]_0$ denotes the commutator 
in the sense of a bounded linear operator on 
$L^{\frac{1}{1-\varepsilon}}(\R^n)$ as before.
All we have to do  is to show that
${\rm I}, {\rm II}\leq C\|f\|_{\mathcal{M}^p_q(dx,w)}$.
For the first term, we notice that 
$[b,T]_0(f\cdot\chi_{2Q})\in L^{\frac{1}{1-\varepsilon}}(\R^n)$, 
particularly,
\begin{equation*}
\lim_{l\to\infty}m_{[b,T]_0(f\cdot\chi_{2Q})}(2^lQ)=0
\end{equation*} 
holds for all $Q\in\mathcal{Q}$ and all medians.
Hence, we may apply Theorem \ref{th-151129-4} 
to obtain that 
\begin{equation*}
{\rm I}
\leq
\left\|
[b,T]_0(f\cdot\chi_{2Q})
\right\|_{\mathcal{M}^p_q(dx,w)}
\leq
C
\left\|
\left\{[b,T]_0(f\cdot\chi_{2Q})
\right\}^\sharp
\right\|_{\mathcal{M}^p_q(dx,w)}.
\end{equation*}
Since we have the pointwise estimate; 
see \cite[Lemma 3.5.5]{grafakos-book}:
\begin{equation*}
\left\{[b,T]_0(f\cdot\chi_{2Q})
\right\}^\sharp(x)
\leq
C
\|b\|_{\rm BMO}
\left(
M^{(\eta)}
[f\cdot\chi_{2Q}](x)
+
M^2[f\cdot\chi_{2Q}](x)
\right)
\end{equation*}
for any $\eta>1$, 
by taking 
$\eta:=\frac{1}{1-\varepsilon}>1$,
it follows that 
\begin{equation*}
{\rm I}
\leq
C
\|b\|_{\rm BMO}
\|f\|_{\mathcal{M}^p_q(dx,w)}.
\end{equation*}
Here, we also used the boundedness of 
$M^{\left(\frac{1}{1-\varepsilon}\right)}$ 
on 
$\mathcal{M}^p_q(dx,w)$.

For the second term, 
we employ Lemma \ref{lm-160114-9}  
to obtain that
\begin{eqnarray*}
{\rm II}
&\leq&
C
\left(
\|b\|_{\rm BMO}
\|f\|_{\mathcal{M}^p_q(dx,w)}
+
\frac{\|f\|_{\mathcal{M}^p_q(dx,w)}}{\Phi_{p,q,w}(Q)}
|Q|^{\frac{1}{p}-\frac{1}{q}}
\left(\int_Q\left|b(x)-b_Q\right|^q w(x)dx\right)^\frac{1}{q}
\right)\\
&=&
C
\left(
\|b\|_{\rm BMO}
\|f\|_{\mathcal{M}^p_q(dx,w)}
+
\|f\|_{\mathcal{M}^p_q(dx,w)}
\left(\frac{1}{w(Q)}
\int_Q\left|b(x)-b_Q\right|^q w(x)dx\right)^\frac{1}{q}
\right).
\end{eqnarray*}
By virtue of Lemma \ref{lm-160115-4}, 
we see that 
${\rm II}\leq C\|b\|_{\rm BMO}\|f\|_{\mathcal{M}^p_q(dx,w)}$.
\end{proof}

Here, we recall a result of the commutators on weighted Morrey spaces
$\mathcal{M}^p_q(w,w)$ of Komori-Shirai type
to compare this result with our results.
\begin{theorem}[\cite{KoSh09}]
Let $1<q\leq p<\infty$, $w\in A_q$.
\begin{enumerate}
\item
Any singular integral operator $T$ is bounded on $\mathcal{M}^p_q(w,w)$.
\item
Assume that $b\in{\rm BMO}$. 
Then the commutator $[b,T]$ generated 
by a singular integral operator $T$ with respect to $b$
is bounded on $\mathcal{M}^p_q(w,w)$.
\end{enumerate}
\end{theorem}

As another application of Theorem \ref{th-151129-3}, 
we consider the boundedness of $M$ 
on the K\"{o}the dual space of $\mathcal{M}^p_q(dx,w)$.
In \cite{Lerner10}, 
Lerner connected the boundedness of maximal operator $M$ 
and the sharp maximal inequality in the framework of Banach function spaces.
Here, we do not give the definitions of the 
Banach function space $X$ and its K\"{o}the dual space $X'$.
We refer \cite{BS88} for these definitions.
\begin{proposition}{\cite[Corollary 4.3]{Lerner10}}\label{lerner}
Let $X$ be a Banach function space 
on which $M$ is bounded.
Then $M$ is bounded on $X'$ 
if and only if there exists $c>0$ such that 
for any $f\in L^0(\R^n)$ such that $\lim_{R\to\infty}f^*(R)=0$, 
\begin{equation*}
\|f\|_{X}\leq
c\|f^\sharp\|_{X}.
\end{equation*} 
\end{proposition}
Unfortunately, 
the Morrey spaces are not Banach function spaces in general; 
see \cite[Example 3.3]{ST15}. 
Nevertheless, 
we still have an analogy of Proposition \ref{lerner} in the framework 
of ``ball Banach function spaces".
We refer \cite{ST15} for the motivation of the notion of ball Banach function spaces and its definitions. 
We may define the K\"{o}the dual space $X'$ of the ball Banach function spaces $X$ 
by the same way as the one of the Banach function space. 
We also note that 
$\mathcal{M}^p_q(dx,w)$ is the ball Banach function spaces as long as 
$M$ is bounded on $\mathcal{M}^p_q(dx,w)$. 

By a similar argument as in \cite{Lerner10}, 
We obtain the Proposition \ref{lerner} for the ball Banach function spaces version.
\begin{proposition}\label{pr-151209-1}
Let $X$ be a ball Banach function space 
on which $M$ is bounded.
Then $M$ is bounded on $X'$ 
if and only if there exists $c>0$ such that 
for any $f\in L^0(\R^n)$ satisfying 
\begin{equation}\label{160121-15}
\lim_{R\to\infty}f^*(R)=0, 
\end{equation}
it holds that 
\begin{equation}\label{151209-1}
\|f\|_{X}\leq
c\|f^\sharp\|_{X}.
\end{equation} 
\end{proposition}

\begin{remark}\label{rm-151209-1}
The assumption
(\ref{160121-15}) in Proposition \ref{pr-151209-1}
implies
(\ref{160121-2}) in Theorem \ref{th-151129-4}.
In fact, if one assumes 
(\ref{160121-15}), 
then we have that
\begin{equation*}
|m_f(2^lQ)|\leq
(f\cdot\chi_{2^lQ})^*(2^{-2}|2^lQ|)
\leq
f^*(2^{-2}|2^lQ|)\to0.
\end{equation*}
As the example
of $f$, where $f:{\mathbb R} \to {\mathbb R}$ is given by,
\[
f=\sum_{j=3}^\infty
\chi_{[j!,j!+1]},
\]
shows that 
the assumption
(\ref{160121-15}) is stronger than
(\ref{160121-2}).
%Another remark%%%%%%%%%%%%%%%%%%%%%%%%%%%%%
\if0
{\color{red}
Does it hold the converse?? 
That is, can we obtain
$\lim_{R\to\infty}f^*(R)=0$
by assuming 
$\lim_{l\to\infty}
m_f(2^lQ)=0$??
$\rightarrow$ This does not hold, namely we have a counter example.}
\fi
%%%%%%%%%%%%%%%%%%%%%%%%%%%%%%%%%%%%%%
\end{remark}

Thanks to Proposition \ref{pr-151209-1} and 
Remark \ref{rm-151209-1}, 
the inequality (\ref{151209-1}) 
for any $f\in L^0(\R^n)$ satisfying 
(\ref{160121-2})
implies
the boundedness of $M$ on $X'$. 
Particularly, we proved (\ref{151209-1}) with 
%$X=\mathcal{M}^p_q(w,w)$ and 
$X=\mathcal{M}^p_q(dx,w)$ 
under the suitable conditions in Corollary \ref{cr-151201-2}.
On the other hand, 
the second author and Tanaka obtained 
the characterization of the K\"{o}the dual space of 
$\mathcal{M}^p_q(dx,w)$ 
in \cite{SaTa14, ST15}. 
%$\mathcal{M}^p_q(w,w)'=\mathcal{H}^{p'}_{q'}(w,w)$ 
%and 
For $\sigma(x)=w(x)^{-\frac{q'}{q}}$, 
we denote all measurable functions $f$
for which norm 
\begin{eqnarray*}
\|f\|_{H^{q',n(1-q/p)}(dx,\sigma)}
&:=&
\left\|w^{-\frac{1}{q}}f\right\|_{H^{q',n(1-q/p)}}
%&=&
%\inf_{b\in \mathfrak{B}_{n(1-q/p)}}
%\left(\int_{\R^n}|f(x)|^{q'}\sigma(x)b(x)^{-\frac{q'}{q}}dx
%\right)^\frac{1}{q'}\\
%&=&
%\inf_{b\in \mathfrak{B}_{n(1-q/p)}}
%\left(\int_{\R^n}|f(x)|^{q'}
%\left(w(x)b(x)\right)^{-\frac{q'}{q}}dx
%\right)^\frac{1}{q'}
\end{eqnarray*}
is finite
by $H^{q',n(1-q/p)}(dx,\sigma)$. 
Then we have that 
$\mathcal{M}^p_q(dx,w)'=H^{q',n(1-q/p)}(dx,\sigma)$ 
when $1<q\leq p<\infty$; 
see \cite{SaTa14, ST15}.
As a result, we obtain
the boundedness of $M$ 
on the K\"{o}the dual space $H^{q',n(1-q/p)}(dx,\sigma)$ 
under the certain conditions:
\begin{theorem}\label{th-160117-1}
Let $1<q\leq p<\infty$ and  assume the boundedness of $M$ on 
$\mathcal{M}^p_q(dx,w)$ and  
the weighted integral condition
(\ref{weight-integral}) for $p,q$ and $w$.
Then $M$ is also bounded on 
$H^{q',n(1-q/p)}(dx,\sigma)$ with $\sigma(x)=w(x)^{-\frac{q'}{q}}$:
\begin{equation*}
\|Mf\|_{H^{q',n(1-q/p)}(dx,\sigma)}
\leq
C
\|f\|_{H^{q',n(1-q/p)}(dx,\sigma)}
\quad
(f\in H^{q',n(1-q/p)}(dx,\sigma)).
\end{equation*}
\end{theorem}

Rerated to Theorem \ref{th-160117-1}, 
Ho \cite{Ho2012} characterized the {\rm BMO} norm in 
the context of general function space including Lebesgue spaces. 
Let us overview the Ho's result. 
Given a Banach function space $X$ 
equipped with a norm $\| \cdot \|_X$, we define the
generalized {\rm BMO} norm
\[
\| b \|_{{\rm BMO}_X}:=
\sup_{Q\in\mathcal{Q}} 
\frac{1}{\| \chi_Q \|_X} \| (b-b_Q)\chi_Q \|_X .
\]
In \cite{Ho2012}, Ho proved that 
if $M$ is bounded on the associate space $X'$, then 
$\| b \|_{{\rm BMO}_X}$ 
and
$\| b \|_{{\rm BMO}}$
are equivalent.
%We remark that 
%Ho's result \cite{Ho2012} has covered
%the authors' one \cite{IzukiRend2010,IzukiSawano2012}. 
%The statements in \cite{IzukiRend2010,IzukiSawano2012} 
%deeply depend on Diening's work \cite{Diening2005} 
%on variable exponent analysis. 
%Meanwhile, Ho's proof is self-contained 
%and obtained as a by-product of the new results about 
%atomic decomposition introduced in \cite{Ho2012}. 
%\end{enumerate}

We can generalize Ho's result
\cite{Ho2012} to the ball Banach function spaces as before.
Hence, by combining Theorem \ref{th-160117-1}
with the idea of Ho, we obtain the following:
\begin{corollary}
Let $1<q\leq p<\infty$ and 
assume the boundedness of $M$ on $\mathcal{M}^p_q(dx,w)$ and 
the weighted integral condition (\ref{weight-integral})
for $p,q$ and $w$.
Define a generalized {\rm BMO} norm with respect to the weighted Morrey space 
of Samko's type
\begin{equation*}
\|b\|_{{\rm BMO}_{\mathcal{M}^p_q(dx,w)}}
:=
\sup_{Q\in\mathcal{Q}} 
\frac{1}{\Phi_{p,q,w}(Q)}
\|(b-b_Q)\chi_Q\|_{\mathcal{M}^p_q(dx,w)}.
\end{equation*}
Then the norm equivalence 
\begin{equation*}
\|b\|_{{\rm BMO}_{\mathcal{M}^p_q(dx,w)}}
\sim
\|b\|_{\rm BMO}
\end{equation*}
holds
for all $b \in {\rm BMO}$.
\end{corollary}


\begin{thebibliography}{999}
\bibitem{AX12}
D. R.~Adams, J.~Xiao, 
Morrey spaces in harmonic analysis, 
Ark. Mat. , {\bf 50} (2012), 
no. 2, 201--230.

\bibitem{BS88}
C. Bennett, R. Sharpley, 
Interpolation of Operators, Academic Press, 1988.

\bibitem{BGST10}
V. I.~Burenkov, V. S.~Guliyev, A.~Serbetci and T. V.~Tararykova, 
Necessary and sufficient condition for the boundedness of genuine 
singular integral operators in local Morrey-Type spaces, Eurasian Math. J. 
{\bf 1}, no. 1, 32--53, 2010.

\bibitem{CF74}
R. R.~Coifman and C.~Fefferman, 
Weighted norm inequalities for maximal functions and singular integrals, Studia Math., 
{\bf 51} 241--250, 1974.

%\bibitem{CHM}
%D.~Cruz-Uribe, SFO, E.~Hern\'andez and J. M.~Martell, 
%Greedy bases in variable Lebesgue spaces, 
%arXiv:1410.1819.

\bibitem{Fujii91}
N.~Fujii, A condition for a two-weight norm inequality for singular integral operators, 
Studia Math., {\bf98} (3) (1991), 175--190.

%%%%%%%%%%%%%%%%%%%%%%%%%%%%%%%%%%%%%%%%%%%%
\if0
\bibitem{Diening2005}
L.~Diening, 
Maximal functions on Musielak-Orlicz spaces and generalized Lebesgue spaces, 
Bull. Sci. Math. {\bf 129} (2005), 657--700.
\fi
%%%%%%%%%%%%%%%%%%%%%%%%%%%%%%%%%%%%%%%%%%%%


%\bibitem{Duoa-Book}
%J.~Duoandikoetxea, 
%Fourier Analysis, 
%Graduate Studies in Math. {\bf 29}, 
%Amer. Math. Soc., Providence, RI, 2001. 

\bibitem{GR85}
J.~Garcia-Cuerva, J. L.~Rubio de Francia, 
Weighted Norm Inequalities and Related Topics, 
North-Holland, Math. Stud., {\bf 116} (1985).

\bibitem{grafakos-book}
L.~Grafakos, Moderen Fourier Analysis, Third Edition, Springer, GTM250.

\bibitem{HNS15}
D. I.~Hakim, E. Nakai and Y. Sawano, 
Generalized fractional maximal operators and vector-valued inequalities on generalized Orlicz-Morrey spaces, 
Revista Matem\'{a}tica Complutense, (2015), 1--32.

\bibitem{Ho2012}
K.-P.~Ho, 
Atomic decomposition of Hardy spaces and characterization of {\rm BMO} 
via Banach function spaces, 
Anal. Math. {\bf 38} (2012), 173--185.

\bibitem{Hytonen11}
T.~Hyt\"{o}nen, Weighted norm inequalities, 
Lecture notes of a course at the University of Helsinki, Winter 2011, 
http://wiki-app.it.helsinki.fi/download/attachments/64424417/weighted.pdf.

\bibitem{Hytonen11-0}
	T. Hyt\"{o}nen, 
	The sharp weighted bound for general Calder\'{o}n-Zygmund operators, 
	Ann. Math. {\bf 175}(3), (2012), 
	1473--1506.	 

%%%%%%%%%%%%%%%%%%%%%%%%%%%%%%%%%%%%%%%%%%%%%
\if0
\bibitem{IzukiRend2010}
M.~Izuki, Boundedness of commutators on Herz spaces with variable exponent, 
Rend. Circ. Mat. Palermo (2) {\bf 59} (2010), 199--213.

\bibitem{IzukiJAA2013}
M.~Izuki, 
Remarks on Muckenhoupt weights with variable exponent, 
J. Anal. Appl. {\bf 11} (2013), 27--42.

\bibitem{IzukiSawano2012}
M.~Izuki and Y.~Sawano, 
Variable Lebesgue norm estimates for {\rm BMO} functtions, 
Czechoslovak Math. J. {\bf 62} (2012), 717--727.

\bibitem{IST2014}
M.~Izuki, Y.~Sawano and Y.~Tsutsui, 
Variable Lebesgue norm estimates for {\rm BMO} functtions. II, 
Anal. Math. {\bf 40} (2014), 215--230.

\fi
%%%%%%%%%%%%%%%%%%%%%%%%%%%%%%%%%%%%%%%%%%%%%

\bibitem{Jawa85}
B.~Jawerth and A.~Torchinsky, Local sharp maximal functions, J. Approx. Theory, {\bf 43} (1965), 231--270.

\bibitem{Jo65}
F.~John, Quasi-isometric mappings, Seminari 1962--1963 di Analisi, Algebra, Geometria e Topologia,
Ist Nazarene Alta Matematica 2 (Edizioni Cremonese, Rome, 1965), 462--473.

%\bibitem{JohnNirenberg}
%J.~John and L.~Nirenberg, 
%On functions of bounded mean oscillation, 
%Comm. Pure Appl. Math. {\bf 14} (1961), 415--426.

%\bibitem{KR} 
%O.~Kov$\acute{\rm{a}}\check{\rm{c}}$ik and J.~R$\acute{\rm{a}}$kosn\'ik, 
%\textit{On spaces $L^{p(x)}$ and $W^{k,p(x)}$}, 
%Czech. Math. J. {\bf 41} (1991), 592--618.

\bibitem{Komori15}
Y.~Komori-Furuya,
Local Good-$\lambda$ Estimate
for the Sharp Maximal Function and Weighted Morrey Space,\\
Journal of Function Spaces
Volume 2015 (2015), Article ID 651825, 4 pages\\
http://dx.doi.org/10.1155/2015/651825.

\bibitem{KoSh09}
Y.~Komori and S.~Shirai,
Weighted Morrey spaces and a singular integral operator,
Math. Nachr. {\bf 282} (2009), no. 2, 219--231.

\bibitem{Lerner98}
A.K.~Lerner, On weighted estimates of non-increasing rearrangements, 
East J. Approx. {\bf 4} (1998), 277--290.

\bibitem{Lerner10}
A.K.~Lerner, 
Some remarks on the Fefferman-Stein inequality, 
J. Anal. Math. {\bf 112} (2010), 329--349.

\bibitem{Lerner13-2}
A.K.~Lerner,
On an estimate of Calder\"{o}n-Zygmund operators by dyadic positive operators,
J. Anal. Math. {\bf 121} (2013), 141--161.

\bibitem{Lerner-Ombrosi10}
A.K.~Lerner, S.~Ombrosi, 
A boundedness criterion for general maximal operators, 
Publ. Mat., {\bf54} (2010) no. 1, 53--71.

\bibitem{MH76}
B.~Muckenhoupt, R. L.~Wheeden, 
Weighted bounded mean oscillation and the Hilbert transform, 
Studia Math., {\bf 54} (1976), 221--237.

\bibitem{Nakai94} 
E. Nakai,
Hardy-Littlewood maximal operator, 
singular integral operators and 
the Riesz potentials on generalized Morrey spaces,
Math. Nachr. {\bf 166} (1994), 95--103.

\bibitem{NS15}
S.~Nakamura, Generalized weighted Morrey spaces and classical operators, 
Math. Nachr., online.

%\bibitem{NS15-2}
%S.~Nakamura, Embedding properties and trace theorem for weighted Morrey spaces, pre-print.

\bibitem{NNS15}
S.~Nakamura, T.~Noi and Y.~Sawano, 
Generalized Morrey spaces and trace operator, 
Sci. China Math.  {\bf 59}  (2016),  no. 2, 281--336.

\bibitem{NS09}
N~Samko, Weighted Hardy and singular operators in Morrey spaces, 
J. Math. Anal. Appl. {\bf 350} (2009), 56--72.

%\bibitem{OV98}
%J.~Orobitg, J.~Verdera, 
%Choquet integrals, Hausdorff content and the Hardy-Littlewood maximal operator, 
%Bull. London Math. Soc., {\bf 30} (1998), no. 2, 145--150.

\bibitem{Perez95}
C.~P\'{e}rez, Endpoint estimates for commutators of singular integral operators, 
J. Funct. Anal., {\bf128} (1995), 163--185.

\if0
\bibitem{Rubio1}
J. L. Rubio de Francia, 
Factorization and extrapolation of weights, 
Bull. Amer. Math. Soc. (N.S.), {\bf 7} (1982), 393--395.

\bibitem{Rubio2}
J. L. Rubio de Francia, 
A new technique in the theory of $A_p$ weights, 
Topics in modern harmonic analysis, Vol. I, II (Turin/Milan, 1982),  571--579, 
Ist. Naz. Alta Mat. Francesco Severi, Rome, 1983. 

\bibitem{Rubio3}
J. L. Rubio de Francia, 
Factorization theory and $A_p$ weights, 
Amer. J. Math. {\bf 106} (1984), 533--547. 
\fi 


%\bibitem{SaSh08}
%Y.~Sawano and S.~Shirai,
%Compact commutators on Morrey spaces with non-doubling measures,
%Georgian Math. J., {\bf 15} (2008), no. 2, 353--376.

\bibitem{SaTa07-2}
Y.~Sawano and H.~Tanaka,
Sharp maximal inequalities and commutators on Morrey spaces
with non-doubling measures,
Taiwanese J. Math. {\bf 11} (2007), no. 4, 1091--1112.

%\bibitem{SaTa10}
% Y.~Sawano and H.~Tanaka,
% Predual spaces of Morrey spaces with non-doubling measures,
% Tokyo J. Math. {\bf 32} (2009), 471--486.

\bibitem{SaTa14}
Y.~Sawano and H.~Tanaka,
Fatou property of predual Morrey spaces with non-doubling measures,
Int. J. Appl. Math. {\bf 27} (2014), no. 3, 283--296.

\bibitem{ST15}
Y.~Sawano, H.~Tanaka,
The Fatou Property of Block spaces,
J. Math. Sci. Univ. Tokyo,
{\bf 22} (2015), 663--683.

\bibitem{St2} 
E.M. Stein,
Harmonic Analysis{\rm :}
Real-Variable Methods, Orthogonality,
and Oscillatory Integrals,
Princeton Univ. Press, 1993.

\bibitem{Ta15}
H.~Tanaka, 
Two-weight norm inequalities on Morrey spaces, 
Ann. Acad. Sci. Fenn. Math. {\bf 40} (2015), 773--791.

\bibitem{St79}
J.O.~Str\"{o}mberg, Bounded mean oscillation with Orlicz norms and duality of Hardy spaces, Indiana
Univ. Math. J. {\bf 28} (1979), 511--544.

%\bibitem{Zo94}
%C.T.~Zorko,
%Morrey space,
%Proc. Amer. Math. Soc. {\bf 98} (1986), 586--592.


\end{thebibliography}
\end{document}